\documentclass[11pt]{amsart}
\usepackage{amsmath,amssymb,enumerate}
\usepackage{graphicx}
\usepackage{mathrsfs}
\usepackage{color}
\usepackage{comment}
\usepackage{amsthm}
\usepackage{hyperref}
\usepackage{tabu}
\usepackage{lipsum}

\newtheorem{theorem}{Theorem}[section]
\newtheorem{proposition}[theorem]{Proposition}
\newtheorem{corollary}[theorem]{Corollary}
\newtheorem{lemma}[theorem]{Lemma}
\theoremstyle{definition}
\newtheorem{definition}[theorem]{Definition}
\theoremstyle{remark}

\newtheorem{remark}[theorem]{Remark}

\newcommand{\CC}{\mathbb{C}}

\newcommand{\FF}{\mathbb{F}}

\newcommand{\RR}{\mathbb{R}}

\newcommand{\cO}{\mathcal O}

\DeclareMathOperator{\diag}{diag}

\DeclareMathOperator{\rank}{rank}

\DeclareMathOperator{\ev}{ev}

\newcommand{\id}{\text{id}}

\newcommand{\pa}[2]{\frac{\partial #1}{\partial #2}}
\newcommand{\paop}[1]{\pa{}{#1}}

\DeclareMathOperator{\Tan}{tan}

\DeclareMathOperator{\sys}{sys}
\DeclareMathOperator{\Hol}{Hol}

\title[Holomorphicity and systole for stable minimal surfaces]{Holomorphicity and sharp systolic estimates for stable minimal surfaces}

\author{Chao Li}
\address{Courant Institute, New York University, 251 Mercer St, New York, NY 10012, USA}
\email{chaoli@nyu.edu}
\author{Yipeng Wang}
\address{Department of Mathematics, Princeton University, Princeton, NJ 08544, USA}
\email{yipeng.wang@princeton.edu}

\begin{document}
\begin{abstract}
    We prove that a complete, parabolic, branched minimal immersion of a Riemann surface $\Sigma$ in $\mathbb{R}^n$ that is stable on all finite-degree covering spaces has image contained in an even dimensional affine subspace and is holomorphic with respect to a constant orthogonal complex structure. Additionally, we prove that such covering stable minimal surfaces do not exist in a manifold with positive isotropic curvature. More quantitatively, for a stable branched minimal immersion $\Sigma$ into a manifold with uniformly positive isotropic curvature, we prove that the systole of $\Sigma$ is bounded logarithmically in terms of its genus. The logarithmic dependence is sharp.
\end{abstract}
\maketitle

\section{Introduction}

An immersion $f: \Sigma \to (M,g)$ is minimal if it is stationary for the area functional, and a minimal immersion is stable if the second variation is nonnegative. Global rigidity theorems and local estimates of stable minimal surfaces are central themes in geometric variational problems. When $(M,g)$ is the Euclidean space $\mathbb{R}^{n+1} $and $\Sigma$ is of codimension one, the rigidity question is called the stable Bernstein problem, and is closely related to the regularity theory of area-minimizing hypersurfaces. Recently, the stable Bernstein problem for minimal hypersurfaces was completely solved: when $n\le 6$, a complete, two-sided, stable minimal immersion $\Sigma^n \to \mathbb{R}^{n+1}$ is flat (cf. \cite{dCP79,FCS80,Pog81,CL24,CL23,CMR24,CLMS26,Maz24,HLW26}); in all higher dimensions, there are non-flat stable minimal hypersurfaces \cite{HS85}. 

In this paper, we study the global rigidity for two-dimensional stable minimal immersions. Throughout, $\Sigma$ is a connected Riemann surface without boundary and $f:\Sigma\to (M,g)$ is a branched conformal minimal immersion. The high co-dimensionality in this scenario allows much more flexible phenomena: for example, every holomorphic curve (with respect to a constant orthogonal complex structure) in an even dimensional subspace of $\mathbb{R}^{n+1}$ is calibrated and hence area-minimizing. A natural question is whether an appropriate converse holds: whether natural stability assumptions imply  holomorphicity of the minimal immersion.

This question has received extensive investigations. An important case to understand is when $(M,g)$ is flat: the Euclidean space or a flat torus $T^n$. Micallef \cite{Micallef1984,Micallef1986} proved that a complete finite-total-curvature stable branched minimal immersion of $\Sigma^2$ in $\mathbb{R}^4$ or $T^4$ is holomorphic. The analogous result holds for genus-zero stable minimal immersions in higher dimensions \cite{Micallef1984}. The question becomes much more delicate when the genus of $\Sigma$ is positive and $n>4$. Indeed, Arezzo--Micallef \cite{ArezzoMicallef2000} and Arezzo--Micallef--Pirola \cite{AMP2002} constructed stable minimal immersion of higher genus surfaces in flat $T^n$ and $\mathbb{R}^n$ that is not holomorphic with respect to any constant complex structure. 

In \cite{Schoen2006}, Schoen proposed seeking holomorphicity under the stronger assumption that $f: \Sigma\to \mathbb{R}^n$ is also stable when lifted to any finite covering space $\hat \Sigma \to \mathbb{R}^n$ - the condition now called \textit{covering stability}. Under this assumption, Fraser--Schoen \cite{FraserSchoen2025} recently established holomorphicity in the genus-one finite-total-curvature case. 

In the first main result, we extend the Fraser--Schoen theorem to surfaces with arbitrary genus.

\begin{theorem}\label{thm:flat}
Let $M$ be either $\mathbb{R}^n$ or a flat torus $\mathbb{R}^n/\Lambda$, and let $f:\Sigma\to M$ be a complete covering stable branched conformal minimal immersion. Assume that $\Sigma$ is parabolic in the Euclidean case and closed in the torus case. Then $f(\Sigma)$ is contained in an even-dimensional affine subspace of $\mathbb{R}^n$, respectively a translate of an even-dimensional flat subtorus of $\mathbb{R}^n/\Lambda$, and $f$ is holomorphic with respect to a constant orthogonal complex structure on this subspace or subtorus.

\end{theorem}

We point out several other recent interesting developments in this direction. Cheng--Karigiannis--Madnick \cite{CKM2025} proved holomorphicity under ordinary stability and a parallel orthogonal complex structure on the normal bundle. Sagman--Thalmaier \cite{SagmanThalmaier2026} treated the plane and some low-genus cases with a different approach. In a related direction, Edelen--Franco Reyna--Minter \cite{EFM26} established holomorphicity of area-minimizing surfaces with density two.

\subsection{Stability  in manifolds with positive isotropic curvature}
We proceed to study rigidity of stable and covering stable minimal surfaces in Riemannian manifolds. An important and well-suited curvature condition is \textit{positive isotropic curvature} (abbreviated as PIC in the sequel). In their seminal paper, Micallef and Moore \cite{MM1988} proved any branched minimal immersion of $S^2$ into a PIC manifold is unstable. This was extended by Fraser \cite{Fraser2003} to rule out sufficiently large stable minimal $T^2$ under the PIC assumption. Recently, Fraser--Schoen \cite{FraserSchoen2025} obtained quantitative systolic estimates for stable minimal $T^2$ in PIC manifolds. We refer the readers to the recent survey article by Fraser \cite{Fraser2026} for these results.

In our second main result, we rule out covering stable minimal surfaces in PIC manifolds.

\begin{theorem}\label{thm:pic}
Let $(M,g)$ be a Riemannian manifold of dimension at least $4$ with positive isotropic curvature. Every conformal branched minimal immersion of a
closed connected oriented surface into $M$ is unstable on some finite unbranched cover.
\end{theorem}

We also seek quantitative estimates on the size of stable minimal surfaces. A recurring theme is that positive curvature forces smallness of stable minimal surfaces. For stable minimal immersions of $T^2$ in PIC manifolds, this was carried out by Fraser--Schoen \cite{FraserSchoen2025}, where they characterized smallness via an explicit systolic inequality, extending the previous result by Fraser \cite{Fraser2003}. We obtain a result in this direction.

\begin{theorem}\label{thm:PICsystole}
    Let $(M,g)$ be a Riemannian manifold of dimension at least $4$, and let $\Sigma$ be a closed connected Riemann surface of genus $g(\Sigma)\geq 1$. Suppose that $(M,g)$ is $\kappa$-PIC and $f:\Sigma \to M$ is a stable conformal minimal immersion, then the systole of $\Sigma$ with respect to the metric $f^*g$ satisfies the inequality
    \[
    \operatorname{sys}(\Sigma)\leq \frac{32}{\sqrt{\kappa}}(1+\log (4g(\Sigma))).
    \]
\end{theorem}

The logarithmic dependence on the genus of $\Sigma$ in Theorem \ref{thm:PICsystole} is sharp: in Section \ref{sec:log-systole} we present an example of a genus $2$ minimal surface that is stable on a sequence of finite covers whose systoles grow logarithmically with their genus. 

\subsection{Proof ideas}
Consider a Riemannian manifold $(M^n, g)$ with connection $D$ and a stable minimal immersion $f: \Sigma \to M$. Let $N$ denote the normal bundle of $\Sigma$ in $M$. The classical idea is to complexify the normal bundle and the index form: indeed, taking $N_\CC = N\otimes \CC$ and choosing isothermal coordinates $z=x+iy$ on $\Sigma$, stability on $\Sigma$ is equivalent to 
\[Q(s,\bar s)= 4\int_\Sigma \left(|D_{\paop{\bar z}}^\perp s|^2 - |D_{\paop z}^{\Tan}  s|^2  - R(\pa{f}{z},s,\pa{f}{\bar z},\bar s) \right) dx\,dy \ge 0, \]
for all sections $s\in \Gamma(N_\CC)$. The fundamental idea of Micallef \cite{Micallef1984}, Micallef--Moore \cite{MM1988}, Fraser \cite{Fraser2003} and Fraser--Schoen \cite{FraserSchoen2025} is to test $Q$ with isotropic sections that are either holomorphic or almost holomorphic.

Our main new ingredient is to twist $N_\CC$ by unitary flat vector bundles. Let $E$ be a unitary flat complex vector bundle of rank $m$. In a local parallel unitary frame $\{e_\alpha\}_{\alpha=1}^m$, write $S=\sum_{\alpha=1}^m s_\alpha\otimes e_\alpha \in\Gamma(N_{\mathbb{C}}\otimes E)$. Define the twisted index form by
\[Q_E (S,\bar S) = 4\int_\Sigma \sum_{\alpha=1}^m \left(|D_{\paop{\bar z}}^\perp s_\alpha|^2 - |D_{\paop{z}}^{\tan} s_\alpha|^2 - R(\pa{f}{z}, s_\alpha, \pa{f}{\bar z}, \bar s_\alpha)  \right) dx\,dy. \]
Suppose first that $E$ is associated with a unitary representation $\rho:\pi_1(\Sigma)\to U(m)$ having finite image. The pullback of $E$ to the finite cover $\pi:\hat \Sigma \to\Sigma$ corresponding to $\ker\rho$ is trivial as a unitary flat bundle. Consequently, in a global parallel unitary frame, $\pi^*S$ becomes an $m$-tuple of normal sections $\widehat{s}_1,\ldots,\widehat{s}_m$ along $f\circ\pi$. Stability of $\hat \Sigma$ therefore gives
\(Q_E(S,\bar S) \ge 0\)
for every smooth compactly supported section $S$ of $N_{\mathbb{C}}\otimes E$.

In fact, the finite-image restriction is unnecessary: surface groups have the property that finite-image unitary representations are dense in the weak-containment topology \cite{LS}. An approximation argument extends the twisted stability inequality to \textit{every} finite-rank unitary flat bundle. On a closed Riemann surface, the Narasimhan--Seshadri correspondence \cite{NS} identifies irreducible unitary representations with stable holomorphic vector bundles of degree zero. This naturally motivates us to construct stable degree zero holomorphic vector bundles $E$ with holomorphic morphisms $\varphi: E\to N_\CC$ (so we may work with $N_\CC\otimes E^*$). We combine this argument with a careful construction $S\in \Gamma(N\otimes E^*)$ from holomorphic isotropic sections to rule out covering stable minimal surfaces in PIC manifolds.

This particular section $S$ can be used to deduce the systolic estimate. The key is to construct an effective covering of $\Sigma$, together with a Lipschitz cutoff function $\eta:\Sigma\to[0,1]$ that vanishes on the boundaries of the pieces. In each parallel frame, the components of $\eta S$ extend by zero to test sections of $N_{\mathbb C}$. The stability inequality, applied to the section $\eta\,S$, gives $\frac{\kappa}{2}\int_\Sigma\eta^2|S|^2\,d\mu\leq\int_\Sigma|d\eta|^2|S|^2\,d\mu$. We therefore want to choose $\eta$ in a way that $\int_\Sigma\eta^2|S|^2\,d\mu$ is comparable to $\int_\Sigma |S|^2\,d\mu$, and with $|d\eta|$ bounded by a multiple of $\frac{1+\log(g(\Sigma))}{\operatorname{sys}(\Sigma)}$.

Analogous problems were studied in graph theoretic problems by Lee--Sidiropoulos \cite[Lemma 2.7]{LeeSidiropoulos2010}. In our differential-geometric setting, we construct $4g(\Sigma)$ minimizing geodesic segments on $\Sigma$ whose complement is a disk. Next, we choose a radius $t\in (0,\frac{\sys(\Sigma)}{4})$ and an ordering of the segments. Varying the radius and ordering gives a family of partitions, and the probabilistic argument of Lee--Sidiropoulos allows us to find a partition whose cutoff satisfies both estimates (See Lemma \ref{lem:localization}). This seems a very general argument that should be applicable to other geometric problems.

The case when $M$ is flat requires a separate argument. For simplicity, we assume $M=T^n$ for the moment. Following Micallef's computations \cite{Micallef1984}, consider a constant vector field $V$, its $(0,1)$-tangential component $V^{(0,1)}$ and its normal component $V^\perp$. The relevant identity (see Lemma \ref{lem:flat-projection-identities} below) is
\[|D^\perp_{\frac{\partial}{\partial\bar z}}V^\perp|^2
 -|D^{\tan}_{\frac{\partial}{\partial z}}V^\perp|^2
 =-\frac{\partial^2}{\partial z\,\partial\bar z}|V^{(0,1)}|^2.\]
Fix points $p,q$, take $z$ the holomorphic coordinate near $p$, and set $V = \pa{f}{z}(p)$. The vector field $z^{-1} V^\perp$ extends smoothly across $p$. If there were a meromorphic function $\sigma$ associated with the divisor $q-p$, the vector field $\sigma V^\perp$ produces a nonpositive second variation. Rigidity then forces that
\[\pa{f}{z}(p) \cdot \pa{f}{z}(q) = 0. \]
Therefore the complex span of all tangent lines is totally isotropic, from which  holomorphicity follows.

This gives an alternative argument when $\Sigma$ has genus zero. When $\Sigma$ has positive genus, the appropriate replacement is a meromorphic section of the degree zero line bundle associated with the divisor $q-p$. This line bundle admits a unitary flat connection, to which we apply the aforementioned twisted index form and still obtain that $\pa{f}{z}(p)\cdot \pa{f}{z}(q) = 0$.

\subsection{AI disclosure}
The authors have worked on problems in this paper for a long time and believed that suitable twisting of the normal bundle is the key. After a series of interactions with a commercial LLM (referred to as `the AI') on covering stable minimal surfaces in $T^n$ in August and early September of 2026, the authors obtained a strategy of proving Theorem \ref{thm:flat}, where two concrete ideas were formed: the flat-projection computational Lemma \ref{lem:flat-projection-identities}, as well as approximating finite dimensional unitary representation by finite-images ones. The authors carried out this strategy to prove the holomorphicity result and extended it to prove the nonexistence of covering stability in PIC (where AI was only used to supply references in algebraic geometry and double check details in the proof of Lemma  \ref{lem:lowering-degree}). 

The most striking aid of AI occured as the authors tried to prove the systolic estimates: after a long series of interactions, the AI found the Lee--Sidiropoulos paper \cite{LeeSidiropoulos2010} in a journal in computer science and is concerned with graph theory (so it would be nearly impossible for the authors to find otherwise). The authors then extended its graph-theoretic proofs in the current differential geometric setting. Finally, the beautiful construction of large covering spaces in Theorem \ref{thm:logarithmic-systole} was obtained in a rough form by the AI - the authors studied and presented a substantially simplified version here. 

This article contains no AI-generated text.

\subsection*{Acknowledgment} The authors are grateful to Rick Schoen for introducing them the question of covering stable minimal surfaces, various stimulating conversations and encouragement over the years. C.L. would also like to thank Valentino Tosatti, Guido De Philippis and Junsheng Zhang for inspiring discussions on related topics, from which he learned a lot. C.L. is supported by NSF grant DMS-2303624 and a Sloan Fellowship.

\section{Preliminaries}

\subsection{Basic notation and definitions}
Let $\Sigma$ be a connected Riemann surface without boundary, let $(M^n,g)$ be a Riemannian manifold, and let $f:\Sigma\to M$ be a branched conformal minimal immersion. In a local conformal coordinate $z=x+iy$, we write $\lambda^2=g(\frac{\partial f}{\partial z},\frac{\partial f}{\partial\bar z})$, so the volume form of $f^*g$ is $d\mu=2\lambda^2\,dx\,dy$. Near a branch point $z=0$, we have $\frac{\partial f}{\partial z}=z^m\,Z$ in a local holomorphic frame, where $m$ is a positive integer and $Z$ is smooth and nonzero at $z=0$. The complexified tangent bundle of $\Sigma$ extends to $z=0$ as the span of $\{Z,\bar Z\}$, and its orthogonal complement defines the complexified normal bundle across the branch point. We write $D^\perp$ and $D^{\tan}$ for the normal and tangential projections of the pullback connection. The operator $D^\perp_{\frac{\partial}{\partial\bar z}}$ determines the holomorphic structure of $N_{\mathbb C}=N\otimes_{\mathbb R}\mathbb C$. The complex bilinear extension of $g$ gives a holomorphic isomorphism $N_{\mathbb C}\simeq N_{\mathbb C}^*$, so $\deg(N_{\mathbb C})=0$ on a closed surface.

Following \cite{Micallef1984} we complexify the index form from the second variation of area. For a compactly supported smooth section $s$ of $N_{\mathbb C}$,
the complexified index form is
\[
Q(s,\bar{s})=4\int_{\Sigma}\Big(|D^{\perp}_{\frac{\partial}{\partial \bar{z}}}s|^2-|D^{\tan}_{\frac{\partial}{\partial{z}}} s|^2-R(\frac{\partial f}{\partial z},s,\frac{\partial f}{\partial \bar z},\bar {s})\Big)\,dx\,dy.
\]
Precisely, setting $s=s_1 + \sqrt{-1}s_2$ for $s_1, s_2$ sections of $N$, we have that $Q(s,\bar s)$ equals the sum of the second variations in the $s_1$ and $s_2$ directions. Note that the integrand extends smoothly across branch points. 
\begin{definition}
    We say that $f:\Sigma\to M$ is stable if $Q(s,\bar s)\geq 0$ for all compactly supported $s\in \Gamma(N_{\mathbb{C}})$. We say $f$ is covering stable if $f$ is stable on every finite cover of $\Sigma$.
\end{definition}

Later we will study the index form on the normal bundle twisted by unitary flat vector bundles. Let $E\to \Sigma$ be an $m$-dimensional unitary flat complex vector bundle with flat unitary connection $\nabla^E$. On a simply connected open set choose a parallel unitary frame $\{e_1,\ldots,e_m\}$ of $E$. For a compactly supported smooth section $S\in \Gamma(N_{\mathbb C}\otimes E)$, write $S=\sum_{\alpha=1}^m s_\alpha\otimes e_\alpha$ in such a frame and define the twisted index form
\[
Q_E(S,\bar S)
=4\int_\Sigma\sum_{\alpha=1}^m\Big(
\big|D^\perp_{\frac{\partial}{\partial\bar z}}s_\alpha\big|^2
-\big|D^{\tan}_{\frac{\partial}{\partial z}} s_\alpha\big|^2-R\Big(\frac{\partial f}{\partial z},s_\alpha,
        \frac{\partial f}{\partial\bar z},\bar s_\alpha\Big)
\Big)\,dx\,dy .
\]
This is well defined globally because two parallel unitary frames differ locally by a constant unitary matrix. Moreover, if $\pi:\widehat\Sigma\to\Sigma$ is a finite covering, then we have $Q_{\pi^*E}(\pi^*S,\overline{\pi^*S})
=\deg(\pi)\,Q_E(S,\bar S)$.

We next recall the following definition for a complete branched minimal immersion $f: \Sigma \to (\mathbb{R}^n,g)$ to be parabolic.
\begin{definition}
    We say that $f: \Sigma\to (\mathbb{R}^n,g)$ is parabolic if there exists a sequence of compactly supported non-negative
smooth functions $\varphi_j:\Sigma\to[0,1]$ such that $\varphi_j=1$ on
each fixed compact subset for all sufficiently large $j$ and
\[
\int_\Sigma |d\varphi_j|^2\,d\mu\longrightarrow 0 .\]
\end{definition}
Finally, when $\Sigma$ is closed, its systole, $\sys(\Sigma)$,  is the length of the shortest noncontractible loop in the metric $f^*g$.

\subsection{The positive isotropic curvature condition}

Given $p\in M$, we define $T_p^{\mathbb{C}}M=T_pM\otimes _{\mathbb{R}}\mathbb{C}$ to be the complexified tangent space to $M$ at $p$. The inner product $g$ extends to a complex bilinear form $g:T_p^{\mathbb{C}}M\times T_p^{\mathbb{C}}M\to \mathbb{C}$, and the Riemann curvature tensor $R$ extends to a complex multilinear form $R:T_p^{\mathbb{C}}M\times T_p^{\mathbb{C}}M\times T_p^{\mathbb{C}}M\times T_p^{\mathbb{C}}M\to \mathbb{C}$. 
\begin{definition}
    Let $\kappa$ be a positive real number. We say $(M,g)$ is $\kappa$-PIC if $$R(\zeta,\eta,\bar\zeta,\bar\eta)\geq \frac{\kappa}{4}\,\big(g(\zeta,\bar{\zeta})\,g(\eta,\bar{\eta})-|g(\zeta,\bar{\eta})|^2\big)$$ for every $p\in M$ and all linearly independent $\zeta,\eta\in T_p^{\mathbb C}M$ satisfying $g(\zeta,\zeta)=g(\zeta,\eta)=g(\eta,\eta)=0$. We say $(M,g)$ has positive isotropic curvature if $R(\zeta,\eta,\bar\zeta,\bar\eta)>0$ for every $p\in M$ and all linearly independent $\zeta,\eta\in T_p^{\mathbb C}M$ satisfying $g(\zeta,\zeta)=g(\zeta,\eta)=g(\eta,\eta)=0$.
\end{definition}

Equivalently, $(M,g)$ is $\kappa$-PIC if for every $p\in M$ and every orthonormal four-frame $\{e_1,e_2,e_3,e_4\}$, we have
\[
R_{1313}+R_{1414}+R_{2323}+R_{2424}-2R_{1234}\geq \kappa,
\]
where $R_{ijk\ell}=R(e_i,e_j,e_k,e_{\ell})$.
\section{Covering Stability and Twisted Index Form}
In this section, we describe a consequence of covering stability in terms of the twisted index form: covering stability implies that $Q_E \ge 0$ for any finite-rank unitary flat complex bundle $E$.
\begin{lemma}\label{lem:finite-image}
    Let $E\to \Sigma$ be an $m$-dimensional unitary flat complex vector bundle with flat unitary
connection $\nabla^E$. Let $\rho:\pi_1(\Sigma)\to U(m)$ be the associated unitary representation. We assume that $\rho(\pi_1(\Sigma))$ is a finite subgroup of $U(m)$. If there exists a compactly supported section $S\in \Gamma(N_{\mathbb{C}}\otimes E)$ such that $Q_E(S,\bar{S})<0$, then $f$ is unstable on some finite cover of $\Sigma$.
\end{lemma}
\begin{proof}
    Let \(\pi:\widehat\Sigma\longrightarrow\Sigma\) be the finite cover corresponding to $\ker\rho$.  Then $\pi^*E$ has a global parallel unitary frame $e_1,\ldots,e_m$.  Next, we write $\pi^*S=\sum_{\alpha=1}^m s_\alpha\otimes e_\alpha$, and it follows that
\[
\deg(\pi)\,Q_E(S,\bar S)
=\sum_{\alpha=1}^m Q(s_\alpha,\bar s_\alpha)<0.
\]
Thus $Q(s_\alpha,\bar s_\alpha)<0$ for some $\alpha\in\{1,\cdots,m\}$.  Hence $f\circ \pi$ is unstable.
\end{proof}
\begin{proposition}\label{prop:index-form-approximation}
Let $\Sigma$ be a connected Riemann surface without boundary. Let
\(
\rho:\pi_1(\Sigma)\to U(V)
\)
be an irreducible finite-dimensional unitary representation, let
$E_\rho\to\Sigma$ be the associated unitary flat bundle, and let
\(
S_{\rho}\in\Gamma(N_{\mathbb C}\otimes E_\rho)
\) be compactly supported.
Then there exist a nonzero vector $v\in V$, a finite set
$K\subset\pi_1(\Sigma)$, and a constant $C_0>0$, depending only on
$\rho$, $S_\rho$, and $v$,  with the following property:

For every finite dimensional unitary representation
\(
\hat{\rho}:\pi_1(\Sigma)\to U(W)
\)
and every $w\in W$, there exists
\(
S_{\hat{\rho}}\in\Gamma(N_{\mathbb C}\otimes E_{\hat{\rho}})
\)
such that
\[
|
Q_{E_\rho}(S_{\rho},\bar S_{\rho})
-
Q_{E_{\hat{\rho}}}(S_{\hat{\rho}},\overline{S}_{\hat{\rho}})
|
\leq
C_0
\max_{\gamma\in K}
|
\langle\rho(\gamma)v,v\rangle
-
\langle\hat{\rho}(\gamma)w,w\rangle
|.
\]
Moreover, $S_{\hat\rho}$ can be chosen so that
$\operatorname{supp}S_{\hat\rho}\subset\operatorname{supp}S_\rho$.
\end{proposition}
\begin{proof}
Let
\(
\pi:\widetilde\Sigma\longrightarrow \Sigma
\)
be the universal covering, and identify
\(
\pi_1(\Sigma)
\)
with the group of deck transformations of $\widetilde\Sigma$. We denote
\(\widetilde N_{\mathbb C}=\pi^*N_{\mathbb C}\) be the pullback of the complexified normal bundle to $\widetilde\Sigma$. 

Next, we write $E_\rho
=
(\widetilde\Sigma\times V)/\pi_1(\Sigma)$ with
\(
\gamma\cdot(x,\xi)
=
(\gamma x,\rho(\gamma)\xi).
\)
Therefore $S_\rho$ lifts to a section
\(
\widetilde S_\rho
\in
\Gamma(\widetilde N_{\mathbb C}\otimes V)
\)
satisfying
\(
\widetilde S_\rho(\gamma x)
=
\rho(\gamma)\,
\widetilde S_\rho(x)
\)
for every $\gamma\in\pi_1(\Sigma)$ and $x\in\widetilde\Sigma$.

We fix a unit vector $v\in V$. Since $\rho$ is irreducible, we have
\(
\operatorname{span}
\{\rho(\gamma)v:\gamma\in\pi_1(\Sigma)\}=V.
\)
Thus, let $m=\dim_{\mathbb C}V$, we may choose
\(
\gamma_1,\ldots,\gamma_m\in\pi_1(\Sigma)
\)
such that
\(
\rho(\gamma_1)v,\ldots,\rho(\gamma_m)v
\)
is a basis of $V$.

We then choose a smooth compactly supported function
\(
\chi:\widetilde\Sigma\to[0,1]
\)
such that
\(
\sum_{\gamma\in\pi_1(\Sigma)}
\chi(\gamma^{-1}x)=1
\) on the preimage of a neighborhood of
$\operatorname{supp}S_\rho$.
We can then find compactly supported sections $s_{\alpha}$ of $\widetilde N_{\mathbb{C}}$ such that
\[
\chi\,\widetilde S_\rho
=
\sum_{\alpha=1}^m
s_\alpha\otimes\rho(\gamma_\alpha)v.
\]
Using the equivariance of $\widetilde S_\rho$, we obtain
\[
\widetilde S_\rho(x)
=
\sum_{\gamma\in\pi_1(\Sigma)}
\sum_{\alpha=1}^m
s_\alpha(\gamma^{-1}x)
\otimes
\rho(\gamma\cdot\gamma_\alpha)v.
\]

Now let
\(
\hat\rho:\pi_1(\Sigma)\to U(W)
\)
be any unitary representation and let $w\in W$. We define
\[
\widetilde S_{\hat\rho}(x)
=
\sum_{\gamma\in\pi_1(\Sigma)}
\sum_{\alpha=1}^m
s_\alpha(\gamma^{-1}x)
\otimes
\hat\rho(\gamma\cdot\gamma_\alpha)w.
\]
Since every $s_\alpha$ has compact support, this sum is locally finite.
Moreover, for every $\zeta\in\pi_1(\Sigma)$,
\begin{align*}
\widetilde S_{\hat\rho}(\zeta x)
&=
\sum_{\gamma\in\pi_1(\Sigma)}\sum_{\alpha=1}^m
s_\alpha(\gamma^{-1}\zeta x)
\otimes
\hat\rho(\gamma\cdot\gamma_\alpha)w\\
&=
\sum_{\eta\in\pi_1(\Sigma)}\sum_{\alpha=1}^m
s_\alpha(\eta^{-1}x)
\otimes
\hat\rho(\zeta\cdot\eta\cdot\gamma_\alpha)w\\
&=
\hat\rho(\zeta)\,\widetilde S_{\hat\rho}(x),
\end{align*}
where in the second equality we set $\gamma=\zeta\cdot\eta$. Hence
$\widetilde S_{\hat\rho}$ descends to a section
\(
S_{\hat\rho}
\in
\Gamma(N_{\mathbb C}\otimes E_{\hat\rho}).
\) Moreover, since $\rho(\gamma_1)v,\ldots,\rho(\gamma_m)v$
form a basis of $V$, the definition of $s_\alpha$
implies
\(
\operatorname{supp}s_\alpha
\subset \pi^{-1}(\operatorname{supp}S_\rho)
\)
for every $\alpha$. Since this set is invariant under deck
transformations, the definition of $\widetilde S_{\hat\rho}$ gives
\(
\operatorname{supp}S_{\hat\rho}
\subset \operatorname{supp}S_\rho.
\)

We now compare the two index forms. Fix a compact fundamental domain
$\Omega\subset\widetilde\Sigma$. Since the $s_\alpha$ have compact support,
there exists a finite set
\(
F\subset\pi_1(\Sigma)
\)
such that, for $x\in\Omega$,
\(
s_\alpha(\gamma^{-1}x)=0
\)
for every $\alpha\in\{1,\cdots,m\}$ whenever $\gamma\notin F$. In particular, for every $\gamma,\eta\in F$ and every $1\leq \alpha,\beta\leq m$, there exists some constant $c_{\gamma,\eta,\alpha,\beta}\in\mathbb C$ depending only on $s_{\alpha}$ and the $f$, such that
\[
Q_{E_\rho}(S_\rho,\overline{S_\rho})
=
\sum_{\gamma,\eta\in F}\sum_{1\leq \alpha,\beta\leq m}
c_{\gamma,\eta,\alpha,\beta}
\left\langle
\rho(\gamma\cdot\gamma_\alpha)v,
\rho(\eta\cdot\gamma_\beta)v
\right\rangle.
\]
Since $\rho$ is unitary, we have
\[
\langle
\rho(\gamma\cdot\gamma_\alpha)v,
\rho(\eta\cdot\gamma_\beta)v
\rangle
=
\langle
\rho(\gamma_\beta^{-1}\cdot\eta^{-1}\cdot\gamma\cdot\gamma_\alpha)v,
v
\rangle.
\]
Therefore, we have
\[
Q_{E_\rho}(S_\rho,\overline{S}_\rho)
=
\sum_{\gamma,\eta\in F}\sum_{1\leq \alpha,\beta\leq m}
c_{\gamma,\eta,\alpha,\beta}
\langle
\rho(\gamma_\beta^{-1}\cdot\eta^{-1}\cdot\gamma\cdot\gamma_\alpha)v,v
\rangle.
\]
Similarly, we have
\[
Q_{E_{\hat\rho}}
(S_{\hat\rho},\overline{S}_{\hat\rho})
=
\sum_{\gamma,\eta\in F}\sum_{1\leq \alpha,\beta\leq m}
c_{\gamma,\eta,\alpha,\beta}
\langle
\hat\rho(\gamma_\beta^{-1}\cdot\eta^{-1}\cdot\gamma\cdot\gamma_\alpha)w,w
\rangle.
\]
We now define
\[
K=
\{
\gamma_\beta^{-1}\cdot\eta^{-1}\cdot\gamma\cdot\gamma_\alpha:
\gamma,\eta\in F,\ 
1\leq\alpha,\beta\leq m
\}
\]
and
\(
C_0=1+
\sum_{\gamma,\eta\in F}\sum_{1\leq \alpha,\beta\leq m}
|c_{\gamma,\eta,\alpha,\beta}|.
\)
Then $K$ is finite, and
\[
\begin{split}
&
|
Q_{E_\rho}(S_\rho,\overline{S}_{\rho})
-
Q_{E_{\hat\rho}}
(S_{\hat\rho},\overline{S}_{\hat\rho})
|
\leq
C_0
\max_{\zeta\in K}
\left|
\langle\rho(\zeta)v,v\rangle
-
\langle\hat\rho(\zeta)w,w\rangle
\right|.
\end{split}
\]
The set $K$ and the constant $C_0$ are independent of $\hat{\rho}$ and $w$. This proves the assertion.
\end{proof}
Next, we state the approximation results regarding representation of surface groups.
\begin{lemma}\label{Lem:FD}
Assume that $\Sigma$ is a closed connected Riemann surface with $g(\Sigma)\geq 1$. Let
\(
\rho:\pi_1(\Sigma)\to U(V)
\)
be an irreducible finite dimensional unitary representation. Given $v\in V$, a finite set
$K\subset\pi_1(\Sigma)$, and $\varepsilon>0$, there exist a unitary finite dimensional representation
\(
\hat{\rho}:\pi_1(\Sigma)\to U(W)
\)
with finite image, and a vector $w\in W$ such that
\[
|
\langle \rho(\gamma)v,v\rangle
-
\langle \hat{\rho}(\gamma)w,w\rangle
|
<\varepsilon
\]
for every $\gamma\in K$.
\end{lemma}

\begin{proof}
This is the consequence of property {\rm FD} for surface groups. See Lubotzky and Y. Shalom \cite{LS}.
\end{proof}
\begin{lemma}
Assume that $\Sigma$ is a closed connected Riemann surface. Let $E\to\Sigma$ be a finite-dimensional
unitary flat complex vector bundle. If there exists
\(
S\in\Gamma(N_{\mathbb C}\otimes E)
\)
such that
\(
Q_E(S,\bar S)<0,
\)
then $f$ is unstable on some finite cover of $\Sigma$.
\end{lemma}
\begin{proof}
Suppose that $Q_E(S,\bar S)<0$ for some finite dimensional unitary flat complex vector bundle $E$ and some $S\in \Gamma(N_{\mathbb{C}}\otimes E)$. If $g(\Sigma)=0$, then $E$ has trivial holonomy and hence contradicting Lemma \ref{lem:finite-image}. We now assume $g(\Sigma)\geq 1$, and let \(\rho:\pi_1(\Sigma)\to U(V)\) be the unitary representation associated with $E$. Since $\rho$ is
finite-dimensional and unitary, we can write
\(
V=V_1\oplus\cdots\oplus V_\ell
\)
as an orthogonal direct sum of irreducible invariant subspaces. Accordingly,
\(
E=E_1\oplus\cdots\oplus E_\ell,
\)
where each $E_j$ is a unitary flat bundle associated with an irreducible
representation
\(
\rho_j:\pi_1(\Sigma)\to U(V_j).
\)

Write
\(
S=S_1+\cdots+S_\ell\) with \(
S_j\in\Gamma(N_{\mathbb C}\otimes E_j).
\)
It then follows that
\[
Q_E(S,\bar S)
=
\sum_{j=1}^\ell Q_{E_j}(S_j,\overline{S}_j).
\]
Since $Q_E(S,\bar S)<0$, there exists some $j\in\{1,\cdots,\ell\}$ such that
\(
Q_{E_j}(S_j,\overline{S}_j)<0.
\)
Thus,
we may assume that $\rho$ is irreducible.

By Proposition \ref{prop:index-form-approximation}, there exist a unit vector $v\in V$,
a finite set
\(
K\subset\pi_1(\Sigma),
\)
and a constant $C_0>0$ such that, for every finite-dimensional unitary
representation
\(
\hat\rho:\pi_1(\Sigma)\to U(W)
\)
and every $w\in W$, there exists
\(
S_{\hat\rho}\in
\Gamma(N_{\mathbb C}\otimes E_{\hat\rho})
\)
satisfying
\[
|
Q_E(S,\bar S)
-
Q_{E_{\hat\rho}}
(S_{\hat\rho},\overline{S}_{\hat\rho})
|
\leq
C_0
\max_{\gamma\in K}
\left|
\langle\rho(\gamma)v,v\rangle
-
\langle\hat\rho(\gamma)w,w\rangle
\right|.
\]
Moreover, Lemma \ref{Lem:FD} implies for any $\varepsilon>0$, there exist a finite-dimensional
unitary representation
\(
\hat\rho:\pi_1(\Sigma)\to U(W)
\)
with finite image and a vector $w\in W$ such that
\[
\max_{\gamma\in K}
|
\langle\rho(\gamma)v,v\rangle
-
\langle\hat\rho(\gamma)w,w\rangle
|
<\varepsilon.
\]
Hence, we have
\[
\begin{aligned}
Q_{E_{\hat\rho}}
(S_{\hat\rho},\overline{S}_{\hat\rho})
&\leq
Q_E(S,\bar S)+C_0\,\varepsilon<0.
\end{aligned}
\]
provided $\varepsilon$ is sufficiently small. Since $\hat\rho$ has finite image, the Lemma \ref{lem:finite-image} implies that $f$ is
unstable on some finite cover of $\Sigma$.
\end{proof}
\begin{corollary}\label{cor:destablization}
Assume that $\Sigma$ is closed and that $f$ is covering stable. Then for
every finite-dimensional unitary flat bundle $E\to\Sigma$ and every
\(
S\in\Gamma(N_{\mathbb C}\otimes E),
\)
we have
\(
Q_E(S,\bar S)\geq0.
\)
\end{corollary}
We now describe analogous results in the non-compact case.
\begin{lemma}\label{lem:non-compact-fundamental-group}
    Let $\Sigma$ be a non-compact Riemann surface without boundary. Let \(\rho:\pi_1(\Sigma)\to U(1)\) be a unitary representation. Given a finite set \(K\subset\pi_1(\Sigma)\) and \(\varepsilon>0\), there exists a character \(\hat\rho:\pi_1(\Sigma)\to U(1)\) with finite image such that
\(
|\rho(\gamma)-\hat\rho(\gamma)|<\varepsilon
\)
for every \(\gamma\in K\).
\end{lemma}
\begin{proof}
Since $\pi_1(\Sigma)$ is free \cite[Chapter I, Theorem 44A]{AS}, we may choose a free basis \(\{\gamma_\alpha\}\) of $\pi_1(\Sigma)$ and write
\(
\rho(\gamma_\alpha)=e^{2\pi i\theta_\alpha}
\)  for some $0\leq\theta_\alpha<1$. Let \(\ell\) be a positive integer such that every element of \(K\) is a product of at most \(\ell\) chosen free generators and their inverses. We choose an integer \(m>2\pi\ell/\varepsilon\) and define
\[
\hat\rho(\gamma_\alpha)=
\begin{cases}
1,
&0\leq\theta_\alpha<\dfrac1m,\\[2mm]
e^{2\pi i k/m},
&\dfrac{k}{m}\leq\theta_\alpha<\dfrac{k+1}{m},
\quad k=1,\ldots,m-1.
\end{cases}
\]
Then $\hat{\rho}$ extends to a representation of $\pi_1(\Sigma)$ whose image is contained in a cyclic group of order $m$. Moreover, for each generator $\gamma_{\alpha}$ we have $|\rho(\gamma_\alpha)-\hat\rho(\gamma_\alpha)|
<\frac{2\pi}{m}$. Hence, for each $\gamma\in K$, we have $|\rho(\gamma)-\hat\rho(\gamma)|<\frac{2\pi\ell}{m}\leq \varepsilon$. This completes the proof. 
\end{proof}

\begin{lemma}\label{lem:non-compact-twist}
    Assume that $\Sigma$ is a noncompact connected Riemann surface and $f$ is covering stable, then for every unitary flat line bundle $L\to \Sigma$ and every compactly supported smooth section $S$ of $N_{\mathbb C}\otimes L$, we have $Q_L(S,\bar{S})\geq 0$.
\end{lemma}
\begin{proof}
    Suppose to the contrary that there exists a unitary flat line bundle $L\to \Sigma$, and and a compactly supported smooth section $S$ of $N_{\mathbb C}\otimes L$ with $Q_L(S,\bar S)<0$. Let \(\rho:\pi_1(\Sigma)\to U(1)\) be the representation defining \(L\). By Lemma \ref{lem:non-compact-fundamental-group} and Proposition \ref{prop:index-form-approximation}, there exist a character \(\hat\rho:\pi_1(\Sigma)\to U(1)\) with finite image and a smooth section \(S_{\hat\rho}\) of $N_{\mathbb C}\otimes E_{\hat{\rho}}$ such that
\(
\operatorname{supp}S_{\hat\rho}\subset\operatorname{supp}S
\)
and \(
Q_{E_{\hat\rho}}(S_{\hat\rho},\overline{S}_{\hat\rho})
<
\frac12 Q_L(S,\bar S)<0.
\)
In particular, \(S_{\hat\rho}\) has compact support. The claim then follows from Lemma \ref{lem:finite-image}.
\end{proof}
\section{Covering Stable Minimal Surface in Flat Space}
In this section, we assume $M$ is either the Euclidean space $\mathbb{R}^n$ or a flat torus $\mathbb{R}^n/\Lambda$ with $n\geq 4$. Let $g$ be a flat metric on $M$. We assume $\Sigma$ is a complete connected Riemann surface, and $f:\Sigma\to M$ is a covering stable branched conformal minimal immersion. We then have that:
\[g\left(\pa{f}{z}, \pa{f}{z}\right)=0,\quad \partial_z\partial_{\bar z}f = 0. \]
We denote by $\lambda^2=g(\frac{\partial f}{\partial z},\frac{\partial f}{\partial \bar z})$. We start prove some computational lemmas following Micallef \cite{Micallef1984}.
\begin{lemma}\label{lem:conformal-factor}
Away from the branch points of $f$, we have
\begin{align*}
D_{\frac{\partial}{\partial z}}^{\tan}\,
\frac{\partial f}{\partial z}
&=
\lambda^{-2}\frac{\partial\lambda^2}{\partial z}
\frac{\partial f}{\partial z},\\
\frac{\partial^2}{\partial z\partial\bar z}\log\lambda^2
&=
\lambda^{-2}
\Big|
D^\perp_{\frac{\partial}{\partial z}}
\frac{\partial f}{\partial z}
\Big|^2.
\end{align*}
\end{lemma}
\begin{proof}
Since $f$ is minimal, we have
\(
D_{\frac{\partial}{\partial\bar z}}
\frac{\partial f}{\partial z}=0.
\)
Hence, we have
\(
\frac{\partial^2\lambda^2}{\partial z\partial\bar z}
=|
D_{\frac{\partial}{\partial z}}
\frac{\partial f}{\partial z}
|^2.
\)
Moreover, conformality gives
\(
g(
D_{\frac{\partial}{\partial z}}
\frac{\partial f}{\partial z},
\frac{\partial f}{\partial z}
)=0,
\)
and therefore we have
\[
D_{\frac{\partial}{\partial z}}^{\tan}
\frac{\partial f}{\partial z}
=
\lambda^{-2}\frac{\partial\lambda^2}{\partial z}
\frac{\partial f}{\partial z}.
\]
Consequently, we have
\[
\Big|
D_{\frac{\partial}{\partial z}}
\frac{\partial f}{\partial z}
\Big|^2
=
\lambda^{-2}
\Big|\frac{\partial\lambda^2}{\partial z}\Big|^2
+
\Big|
D^\perp_{\frac{\partial}{\partial z}}
\frac{\partial f}{\partial z}
\Big|^2.
\]
From this, the assertion follows.
\end{proof}
In the rest of this section, for every constant vector $V\in\mathbb C^n$, we write $V=V^{\tan}+V^{\perp}$ for the tangential and normal decomposition. In particular, $V^{\perp}\in\Gamma(N_{\mathbb{C}})$. Moreover, we define $V^{\tan}=V^{(1,0)}+V^{(0,1)}$ by
\begin{align*}
    V^{(1,0)}&=\lambda^{-2}g(V,\frac{\partial f}{\partial \bar z})\,\frac{\partial f}{\partial z},\\
    V^{(0,1)}&=\lambda^{-2}g(V,\frac{\partial f}{\partial  z})\,\frac{\partial f}{\partial \bar z}.
\end{align*}
Near a branch point $z=0$, we write $\frac{\partial f}{\partial z}=z^m\,Z$ with $Z(0)\neq 0$. Then $V^{(1,0)}$ and $V(0,1)$ extend smoothly across $z=0$ by $V^{(1,0)}=g(Z,\bar{Z})^{-1}g(V,\bar Z)\,Z$ and $V^{(0,1)}=g(Z,\bar{Z})^{-1}g(V, Z)\,\bar Z$.

The next lemma is the key computation.

\begin{lemma}\label{lem:flat-projection-identities}
For every constant $V\in\mathbb C^n$, we have
\begin{align*}
g(D^\perp_{\frac{\partial}{\partial\bar z}}V^\perp,
          \overline{V^\perp})
 &=-\frac{\partial}{\partial\bar z}|V^{(0,1)}|^2,\\
|D^\perp_{\frac{\partial}{\partial\bar z}}V^\perp|^2
 -|D^{\tan}_{\frac{\partial}{\partial z}}V^\perp|^2
 &=-\frac{\partial^2}{\partial z\,\partial\bar z}|V^{(0,1)}|^2 .
\end{align*}
\end{lemma}

\begin{proof}
Since $f$ is a conformal minimal immersion, we have
\[
D_{\frac{\partial}{\partial\bar z}}\frac{\partial f}{\partial z}
=D_{\frac{\partial}{\partial z}}\frac{\partial f}{\partial\bar z}=0,
\]
and
\[
\frac{\partial}{\partial\bar z}g(V,\frac{\partial f}{\partial z})
=\frac{\partial}{\partial z}g(V,\frac{\partial f}{\partial\bar z})=0.
\]
Using Lemma \ref{lem:conformal-factor}, we can differentiate the identity $V^\perp=V-V^{(1,0)}-V^{(0,1)}$ and obtain
\begin{align*}
&D^\perp_{\frac{\partial}{\partial\bar z}}V^\perp\\
 &=-\lambda^{-2}g(V,\frac{\partial f}{\partial\bar z})
       D^\perp_{\frac{\partial}{\partial\bar z}}\frac{\partial f}{\partial z}-\lambda^{-2}g(V,\frac{\partial f}{\partial z})
       D^\perp_{\frac{\partial}{\partial\bar z}}\frac{\partial f}{\partial\bar z}\\
 &=-\lambda^{-2}g(V,\frac{\partial f}{\partial z})
       D^\perp_{\frac{\partial}{\partial\bar z}}\frac{\partial f}{\partial\bar z}.
\end{align*}
Moreover,
\[
\frac{\partial}{\partial z}
  \Big(\lambda^{-2}\,g(V,\frac{\partial f}{\partial\bar z})\Big)
=-\lambda^{-4}\frac{\partial\lambda^2}{\partial z}
       g(V,\frac{\partial f}{\partial\bar z}),
\]
Combining this with Lemma \ref{lem:conformal-factor}, we obtain
\begin{align*}
&D^{\tan}_{\frac{\partial}{\partial z}}V^\perp\\
&=-\frac{\partial}{\partial z}
       \Big(\lambda^{-2}g(V,\frac{\partial f}{\partial\bar z})\Big)\,
                              \frac{\partial f}{\partial z}-\lambda^{-4}g(V,\frac{\partial f}{\partial\bar z})\,
       \frac{\partial\lambda^2}{\partial z}\,\frac{\partial f}{\partial z}\\
 &\qquad -\frac{\partial}{\partial z}
       \Big(\lambda^{-2}g(V,\frac{\partial f}{\partial z})\Big)\,\frac{\partial f}{\partial\bar z}\\
&=-\frac{\partial}{\partial z}
       \Big(\lambda^{-2}g(V,\frac{\partial f}{\partial z})\Big)\,
                              \frac{\partial f}{\partial\bar z}.
\end{align*}
Using Lemma \ref{lem:conformal-factor} again, we have
\begin{align*}
\lambda^{-2}D^\perp_{\frac{\partial}{\partial\bar z}}
                          \frac{\partial f}{\partial\bar z}
 =\lambda^{-2}D_{\frac{\partial}{\partial\bar z}}
                          \frac{\partial f}{\partial\bar z}
   -\lambda^{-4}\frac{\partial\lambda^2}{\partial\bar z}
                          \frac{\partial f}{\partial\bar z}=D_{\frac{\partial}{\partial\bar z}}
              \Big(\lambda^{-2}\frac{\partial f}{\partial\bar z}\Big),
\end{align*}
This implies
\begin{align*}
&g(D^\perp_{\frac{\partial}{\partial\bar z}}V^\perp,
           \overline{V^\perp})\\
&=-\lambda^{-2}\,g(V,\frac{\partial f}{\partial z})\,
       g(\bar V,D^\perp_{\frac{\partial}{\partial\bar z}}
                         \frac{\partial f}{\partial\bar z})\\
&=-g(V,\frac{\partial f}{\partial z})\,
       g(\bar V,D_{\frac{\partial}{\partial\bar z}}
                  \Big(\lambda^{-2}\frac{\partial f}{\partial\bar z}\Big))\\
&=-g(V,\frac{\partial f}{\partial z})\,
       \frac{\partial}{\partial\bar z}
           \Big(\lambda^{-2}g(\bar V,\frac{\partial f}{\partial\bar z})\Big)\\
&=-\frac{\partial}{\partial\bar z}
       \Big(\lambda^{-2}\,g(V,\frac{\partial f}{\partial z})
                         g(\bar V,\frac{\partial f}{\partial\bar z})\Big)\\
&=-\frac{\partial}{\partial\bar z}|V^{(0,1)}|^2.
\end{align*}
Next, we compute
\begin{align*}
\frac{\partial}{\partial z}|V^{(0,1)}|^2
 &=g(\bar V,\frac{\partial f}{\partial\bar z})
       \frac{\partial}{\partial z}
          \Big(\lambda^{-2}g(V,\frac{\partial f}{\partial z})\Big)\\
 &=\lambda^2\frac{\partial}{\partial z}
          \Big(\lambda^{-2}g(V,\frac{\partial f}{\partial z})\Big)
          \lambda^{-2}g(\bar V,\frac{\partial f}{\partial\bar z}).
\end{align*}
Moreover, we have
\begin{align*}
\lambda^2\frac{\partial}{\partial z}
       \Big(\lambda^{-2}g(V,\frac{\partial f}{\partial z})\Big)
 &=\frac{\partial}{\partial z}g(V,\frac{\partial f}{\partial z})
   -\frac{\partial\log\lambda^2}{\partial z}
           g(V,\frac{\partial f}{\partial z}),\\
\frac{\partial}{\partial\bar z}
 \Big(\lambda^2\frac{\partial}{\partial z}
       \Big(\lambda^{-2}g(V,\frac{\partial f}{\partial z})\Big)\Big)
 &=-g(V,\frac{\partial f}{\partial z})
        \,\frac{\partial^2\log\lambda^2}{\partial z\,\partial\bar z}.
\end{align*}
Combining these together and using Lemma \ref{lem:conformal-factor} again, we conclude that
\begin{align*}
&\frac{\partial^2}{\partial z\,\partial\bar z}|V^{(0,1)}|^2\\
&=\lambda^2\,
    \frac{\partial}{\partial z}
       \Big(\lambda^{-2}g(V,\frac{\partial f}{\partial z})\Big)\,
    \frac{\partial}{\partial\bar z}
       \Big(\lambda^{-2}g(\bar V,\frac{\partial f}{\partial\bar z})\Big)-|V^{(0,1)}|^2\,
       \frac{\partial^2\log\lambda^2}{\partial z\,\partial\bar z}\\
&=\lambda^2\,\Big|
    \frac{\partial}{\partial z}
       \Big(\lambda^{-2}g(V,\frac{\partial f}{\partial z})\Big)
       \Big|^2-\lambda^{-4}\Big|g(V,\frac{\partial f}{\partial z})\Big|^2\cdot
       \big|D^\perp_{\frac{\partial}{\partial z}}
                        \frac{\partial f}{\partial z}\big|^2\\
&=|D^{\tan}_{\frac{\partial}{\partial z}}V^\perp|^2
    -|D^\perp_{\frac{\partial}{\partial\bar z}}V^\perp|^2.
\end{align*}
This completes the proof.
\end{proof}
Fix distinct points $p,q\in\Sigma$ away from the branch points, choose a conformal coordinate $z$ centered at $p$, and set $V=\frac{\partial f}{\partial z}(p)$. Let $L=\mathcal O(q-p)$ and let $\sigma$ be its meromorphic section with a simple pole at $p$ and a simple zero at $q$. If $\Sigma$ is closed, the degree-zero line bundle $L$ admits a flat Hermitian metric.
\begin{lemma}
    If $\Sigma$ is non-compact, then there exists a unitary flat metric on $L$ such that $|\sigma|$ is bounded outside a compact subset of $\Sigma$.
\end{lemma}
\begin{proof}
    Choose a conformal coordinate \(w\) centered at \(q\). By \cite[Chapter III, \S3, Theorem 9A]{AS}, there exists a real harmonic function \(G\) on \(\Sigma\setminus\{p,q\}\), bounded outside a compact subset, such that $G+\log|z|$, $G-\log|w|$
extend harmonically across \(p\) and \(q\), respectively. We then define a Hermitian metric on $L|_{\Sigma\setminus \{p,q\}}$ such that $|\sigma|=e^G$. Since \(z\sigma\) and \(w^{-1}\sigma\) are nonvanishing holomorphic frames near \(p\) and \(q\), respectively. One also has 
\(
\log|z\sigma|=G+\log|z|\) and \(
\log|w^{-1}\sigma|=G-\log|w|.
\)
Therefore, the Hermitian metric extends smoothly across $p$ and $q$. Clearly, the associated Chern connection is flat, and $|\sigma|$ is bounded outside a compact subset of $\Sigma$. 
\end{proof}
\begin{lemma}\label{lem:vanishing-order}
$z^{-1}V^\perp$ extends smoothly across $p$ as a section
of $N_{\mathbb C}$. Moreover, we have
$|V^{(0,1)}|^2=O(|z|^4)$ and $d|V^{(0,1)}|^2=O(|z|^3)$ near $p$.
\end{lemma}

\begin{proof}
Since $f$ is minimal, we know that
$\frac{\partial f}{\partial z}$ is a holomorphic
$\mathbb C^n$-valued function. In particular, since
$V=\frac{\partial f}{\partial z}(p)$, we know that
\(
z^{-1}(V-\frac{\partial f}{\partial z})
\)
extends holomorphically across $p$ as a $\mathbb C^n$-valued
function.

On the other hand, the orthogonal projection from
$\Sigma\times\mathbb C^n$ onto $N_{\mathbb C}$ is a smooth
complex-linear map. Since
$\frac{\partial f}{\partial z}$ is tangent to $\Sigma$,
we have
\(
z^{-1}V^\perp
=(z^{-1}(V-\frac{\partial f}{\partial z}))^\perp
\)
extends
smoothly across $p$ as a section of $N_{\mathbb C}$.

Moreover, conformality gives
\(
g\big(V,\frac{\partial f}{\partial z}\big)
=-\frac12g\big(
V-\frac{\partial f}{\partial z},
V-\frac{\partial f}{\partial z}\big).
\)
Since $V-\frac{\partial f}{\partial z}$ vanishes to
order at least one at $p$, we conclude that
$g(V,\frac{\partial f}{\partial z})$ vanishes to order
at least two at $p$. Thus,
\(
|V^{(0,1)}|^2
=\lambda^{-2}
|g(V,\frac{\partial f}{\partial z})|^2
=O(|z|^4)
\) around $p$. Differentiating again gives $d|V^{(0,1)}|^2=O(|z|^3)$ around $p$.
This completes the proof.
\end{proof}
In particular, Lemma \ref{lem:vanishing-order} implies $S=V^{\perp}\otimes \sigma$ extends smoothly across $p$ to a section of $N_{\mathbb{C}}\otimes L$.
\begin{lemma}\label{lem:flat-cutoff}
For every real smooth compactly supported function $\varphi$ on $\Sigma$,
\begin{align*}
&Q_L(\varphi S,\varphi\bar S)\\
=&\int_\Sigma |d\varphi|^2\,|\sigma|^2\,
          \bigl(|V^\perp|^2+|V^{(0,1)}|^2\bigr)\,d\mu\\
          &-4\int_\Sigma |V^{(0,1)}|^2\cdot
      |\nabla^L_{\frac{\partial}{\partial z}}(\varphi\sigma)|^2
      \,dx\,dy .
\end{align*}
\end{lemma}
\begin{proof}
On $\Sigma\setminus\{p\}$, we have
\(
\nabla^L_{\frac{\partial}{\partial\bar z}}\sigma=0\) and \(
\nabla^L_{\frac{\partial}{\partial\bar z}}
 \nabla^L_{\frac{\partial}{\partial z}}\sigma
=\nabla^L_{\frac{\partial}{\partial z}}
 \nabla^L_{\frac{\partial}{\partial\bar z}}\sigma=0.
\)
This implies
\[
\frac{\partial^2|\sigma|^2}{\partial z\,\partial\bar z}
=|\nabla^L_{\frac{\partial}{\partial z}}\sigma|^2.
\]
Since $S=V^\perp\otimes\sigma$, we have
\begin{align*}
D^\perp_{\frac{\partial}{\partial\bar z}}(\varphi S)
 &=\frac{\partial\varphi}{\partial\bar z}\,V^\perp\otimes\sigma
   +\varphi \,D^\perp_{\frac{\partial}{\partial\bar z}}V^\perp\otimes\sigma,\\
D^{\tan}_{\frac{\partial}{\partial z}}(\varphi S)
 &=\varphi \,D^{\tan}_{\frac{\partial}{\partial z}}V^\perp\otimes\sigma .
\end{align*}
By Lemma \ref{lem:flat-projection-identities},
\begin{align*}
&|D^\perp_{\frac{\partial}{\partial\bar z}}(\varphi S)|^2
 -|D^{\tan}_{\frac{\partial}{\partial z}}(\varphi S)|^2\\
&=|\sigma|^2\,\Big|\frac{\partial\varphi}{\partial\bar z}\Big|^2\,|V^\perp|^2
 +\varphi^2|\sigma|^2
    (|D^\perp_{\frac{\partial}{\partial\bar z}}V^\perp|^2
          -|D^{\tan}_{\frac{\partial}{\partial z}}V^\perp|^2)\\
&\quad+2\varphi\,|\sigma|^2\,\operatorname{Re}
    \Big(\frac{\partial\varphi}{\partial z}\,
          g(D^\perp_{\frac{\partial}{\partial\bar z}}V^\perp,
                       \overline{V^\perp})\Big)\\
&=|\sigma|^2\,\Big|\frac{\partial\varphi}{\partial\bar z}\Big|^2\,|V^\perp|^2
 -\varphi^2\,|\sigma|^2\,
       \frac{\partial^2|V^{(0,1)}|^2}{\partial z\,\partial\bar z}\\
&\quad-2\varphi\,|\sigma|^2\,\operatorname{Re}
    \Big(\frac{\partial\varphi}{\partial z}
          \frac{\partial|V^{(0,1)}|^2}{\partial\bar z}\Big)\\
&=|\sigma|^2\,\Big|\frac{\partial\varphi}{\partial\bar z}\Big|^2\,|V^\perp|^2
 -|\sigma|^2\,\operatorname{Re}\frac{\partial}{\partial z}
       \Big(\varphi^2\frac{\partial|V^{(0,1)}|^2}{\partial\bar z}\Big).
\end{align*}
On the other hand, we have
\begin{align*}
&|\nabla^L_{\frac{\partial}{\partial z}}(\varphi\sigma)|^2\\
&=\Big|\frac{\partial\varphi}{\partial z}\sigma
          +\varphi\nabla^L_{\frac{\partial}{\partial z}}\sigma\Big|^2\\
&=|\sigma|^2\,\Big|\frac{\partial\varphi}{\partial\bar z}\Big|^2
 +\varphi^2\frac{\partial^2|\sigma|^2}{\partial z\,\partial\bar z}
 +2\varphi\,\operatorname{Re}
    \Big(\frac{\partial\varphi}{\partial\bar z}
          \frac{\partial|\sigma|^2}{\partial z}\Big)\\
&=|\sigma|^2\,\Big|\frac{\partial\varphi}{\partial\bar z}\Big|^2
 +\operatorname{Re}\frac{\partial}{\partial\bar z}
       \Big(\varphi^2\frac{\partial|\sigma|^2}{\partial z}\Big).
\end{align*}
Around $p$, Lemma \ref{lem:vanishing-order} implies for $r$ sufficiently small, we have
\begin{align*}
&\Big|\int_{|z|=r}\varphi^2\,|\sigma|^2\,
       \frac{\partial|V^{(0,1)}|^2}{\partial\bar z}\,d\bar z\Big|\leq 2\pi r\sup_{|z|=r}
       \Big|\varphi^2\,|\sigma|^2\,
             \frac{\partial|V^{(0,1)}|^2}{\partial\bar z}\Big|
 =O(r^2),\\
&\Big|\int_{|z|=r}\varphi^2\,|V^{(0,1)}|^2\,
       \frac{\partial|\sigma|^2}{\partial z}\,dz\Big|\leq 2\pi r\sup_{|z|=r}
       \Big|\varphi^2\,|V^{(0,1)}|^2\,
             \frac{\partial|\sigma|^2}{\partial z}\Big|
 =O(r^2).
\end{align*}
Furthermore, around $p$ we have $|\sigma|^2(|V^\perp|^2+|V^{(0,1)}|^2)=O(1)$ and $|V^{(0,1)}|^2\cdot
 |\nabla^L_{\frac{\partial}{\partial z}}(\varphi\sigma)|^2=O(1)$. The integrands extend across the branch points away from $p$. Integrating over $\Sigma\setminus \{|z|<r\}$ and taking $r\to0$ gives
\begin{align*}
&\frac14\,Q_L(\varphi S,\varphi\bar S)\\
&=\int_\Sigma|\sigma|^2\,
       \Big|\frac{\partial\varphi}{\partial\bar z}\Big|^2\,|V^\perp|^2\,dx\,dy-\operatorname{Re}\int_\Sigma|\sigma|^2\,
       \frac{\partial}{\partial z}
          \Big(\varphi^2\frac{\partial|V^{(0,1)}|^2}{\partial\bar z}\Big)
                \,dx\,dy\\
&=\int_\Sigma|\sigma|^2\,
       \Big|\frac{\partial\varphi}{\partial\bar z}\Big|^2\,|V^\perp|^2\,dx\,dy+\operatorname{Re}\int_\Sigma\varphi^2\,
       \frac{\partial|\sigma|^2}{\partial z}\,
       \frac{\partial|V^{(0,1)}|^2}{\partial\bar z}\,dx\,dy\\
&=\int_\Sigma|\sigma|^2\,
       \Big|\frac{\partial\varphi}{\partial\bar z}\Big|^2\,|V^\perp|^2\,dx\,dy-\int_\Sigma|V^{(0,1)}|^2\operatorname{Re}
       \frac{\partial}{\partial\bar z}
          \Big(\varphi^2\frac{\partial|\sigma|^2}{\partial z}\Big)\,dx\,dy\\
&=\int_\Sigma|\sigma|^2\,
       \Big|\frac{\partial\varphi}{\partial\bar z}\Big|^2\,
          \bigl(|V^\perp|^2+|V^{(0,1)}|^2\bigr)\,dx\,dy-\int_\Sigma|V^{(0,1)}|^2\cdot
       |\nabla^L_{\frac{\partial}{\partial z}}(\varphi\sigma)|^2
                \,dx\,dy.
\end{align*}
Finally,
\[
4\Big|\frac{\partial\varphi}{\partial\bar z}\Big|^2\,dx\,dy
=\Big(\Big|\frac{\partial\varphi}{\partial x}\Big|^2
       +\Big|\frac{\partial\varphi}{\partial y}\Big|^2\Big)\,dx\,dy
=|d\varphi|^2\,d\mu .
\]
This completes the proof of the Lemma.
\end{proof}

\begin{lemma}\label{lem:isotropic}
If $\Sigma$ is closed or parabolic, then
$g(\frac{\partial f}{\partial z}(p),\frac{\partial f}{\partial z})=0$
on $\Sigma$ for every $p\in\Sigma$.
\end{lemma}

\begin{proof}
Since $\frac{\partial f}{\partial z}=0$ at every branch point, it suffices the prove the claim away from the branch points. If $\Sigma$ is closed, Corollary \ref{cor:destablization} and
Lemma \ref{lem:flat-cutoff} with $\varphi=1$ give
\[
0\leq Q_L(S,\bar S)
=-4\int_\Sigma |V^{(0,1)}|^2
 \Big|\nabla^L_{\frac{\partial}{\partial z}}\sigma\Big|^2\,dx\,dy
\leq0.
\]
If $\Sigma$ is noncompact and parabolic, we let $\varphi_j:\Sigma\to[0,1]$ be a sequence of non-negative smooth compactly supported functions such that $\varphi_j=1$ on each fixed compact subset for all sufficiently large $j$ and
\(
\int_\Sigma |d\varphi_j|^2\,d\mu\longrightarrow 0 \). Since
$\operatorname{supp}(\varphi_jS)\subset\operatorname{supp}\varphi_j$,
Lemma \ref{lem:non-compact-twist} implies
$Q_L(\varphi_jS,\varphi_j\bar S)\geq0$ for each $j$.

By Lemma \ref{lem:vanishing-order}, around $p$ we have
\[
|\sigma|^2\bigl(|V^\perp|^2+|V^{(0,1)}|^2\bigr)
=O(|z|^{-2})\bigl(O(|z|^2)+O(|z|^4)\bigr)=O(1).
\]
Outside a compact set, $|\sigma|$ is bounded and
\[
|\sigma|^2\,(|V^\perp|^2+|V^{(0,1)}|^2)
=|\sigma|^2(|V|^2-|V^{(1,0)}|^2)
\leq|\sigma|^2|V|^2 .
\]
Thus
$\sup_\Sigma|\sigma|^2\,(|V^\perp|^2+|V^{(0,1)}|^2)<\infty$.
Lemma \ref{lem:flat-cutoff} yields
\begin{align*}
0&\leq
4\int_\Sigma |V^{(0,1)}|^2
\cdot|\nabla^L_{\frac{\partial}{\partial z}}(\varphi_j\sigma)|^2\,dx\,dy\\
&\leq\int_\Sigma |d\varphi_j|^2\,|\sigma|^2\cdot
                    \bigl(|V^\perp|^2+|V^{(0,1)}|^2\bigr)\,d\mu\\
&\leq\sup_\Sigma(|\sigma|^2
                    \bigl(|V^\perp|^2+|V^{(0,1)}|^2))\,
                    \int_\Sigma|d\varphi_j|^2\,d\mu
\end{align*}
for each $j\geq 1$. Since $\varphi_j=1$ on every fixed compact set for all sufficiently
large $j$, we conclude that
\(
\nabla^L_{\frac{\partial}{\partial z}}(\varphi_j\sigma)
\longrightarrow\nabla^L_{\frac{\partial}{\partial z}}\sigma\) on  $\Sigma\setminus\{p\}.$
Fatou's lemma gives
\begin{align*}
0&\leq4\int_\Sigma |V^{(0,1)}|^2
       \Big|\nabla^L_{\frac{\partial}{\partial z}}\sigma\Big|^2\,dx\,dy\\
&\leq\liminf_{j\to\infty}4\int_\Sigma |V^{(0,1)}|^2
       \Big|\nabla^L_{\frac{\partial}{\partial z}}(\varphi_j\sigma)\Big|^2
       \,dx\,dy=0 .
\end{align*}
To summarize, if $\Sigma$ is closed or parabolic, we have
\(
|V^{(0,1)}|^2
 \cdot|\nabla^L_{\frac{\partial}{\partial z}}\sigma|^2=0\)
 on $\Sigma\setminus\{p\}$. Choose a conformal coordinate $w$ centered at $q$, since $\nabla^L_{\frac{\partial}{\partial w}}\sigma$ is non-zero in a neighborhood of $q$, we must have $|V^{(0,1)}|^2=0$ near $q$. From this, we conclude that
$g(V,\frac{\partial f}{\partial w})=0$ in a neighborhood of $q$. Since the local expressions
$g(V,\frac{\partial f}{\partial z})\,dz$ define a holomorphic one-form, we conclude that
$g(V,\frac{\partial f}{\partial z})=0$ at each point of $\Sigma$. This completes the proof.
\end{proof}

\begin{proposition}\label{prop:flat-rigidity}
If $\Sigma$ is closed or parabolic, there is an even-dimensional real
subspace $V\subset\mathbb R^n$ with a constant orthogonal complex
structure $J$ such that $df$ takes values in $V$ and
$df\circ j=J\circ df$, where $j$ is the complex structure on $\Sigma$.
\end{proposition}

\begin{proof}
Let
\(
W=\operatorname{span}_{\mathbb C}
       \{\frac{\partial f}{\partial z}(p):p\in\Sigma\}
       \subset\mathbb C^n.
\)
Lemma \ref{lem:isotropic} gives $g(v,w)=0$ for $v,w\in W$. If
$v\in W\cap\overline W$, then $\bar v\in W$ and
\(
|v|^2=g(v,\bar v)=0.
\)
Thus, we have $W\cap\overline W=\{0\}$. We define $V=\{v+\bar v:v\in W\}\subset\mathbb R^n$ so that $V\otimes_{\mathbb R}\mathbb C=W\oplus\overline W$.
 Moreover, we define $J:V\to V$ by
\(
J(v+\bar v)=iv-i\bar v.
\)
It follows that for $v,w\in W$, we have $J^2(v+\bar v)=-v-\bar v$ and
\begin{align*}
g(J(v+\bar v),J(w+\bar w))
 &=g(iv-i\bar v,iw-i\bar w)\\
 &=g(v,\bar w)+g(\bar v,w)\\
 &=g(v+\bar v,w+\bar w).
\end{align*}
Therefore, $J$ is orthogonal and
$\dim_{\mathbb R}V=2\dim_{\mathbb C}W$. In a conformal coordinate, we have
\[
\frac{\partial f}{\partial x}
=\frac{\partial f}{\partial z}+\frac{\partial f}{\partial\bar z},
\qquad
\frac{\partial f}{\partial y}
=i\frac{\partial f}{\partial z}-i\frac{\partial f}{\partial\bar z},
\]
This implies $J\frac{\partial f}{\partial x}=\frac{\partial f}{\partial y}$, $J\frac{\partial f}{\partial y}=-\frac{\partial f}{\partial x}$ and $J\circ df=df\circ j$.
\end{proof}

\begin{proof}[Proof of Theorem \ref{thm:flat} when $M=\RR^n/\Lambda$]
Take $V,J$ from Proposition \ref{prop:flat-rigidity} and let
$\widetilde f:\widetilde\Sigma\to\mathbb R^n$ be a lift of $f$.
Since $d\widetilde f$ takes values in $V$, for every $x,y\in\widetilde{\Sigma}$, we have
$\widetilde f(y)-\widetilde f(x)
=\int_x^y d\widetilde f\in V$. Hence, for every $\gamma\in\pi_1(\Sigma)$, the period $\alpha_{\gamma}=\int_\gamma df$ satisfies $\alpha_{\gamma}\in V\cap \Lambda$.

Suppose $v\in V$ and $\langle v,a_\gamma\rangle=0$ for every $\gamma$.
Then it follows that for every $x\in \tilde{\Sigma}$, we have
\[
\langle v,\widetilde f(\gamma x)\rangle
-\langle v,\widetilde f(x)\rangle
=\langle v,a_\gamma\rangle=0.
\]
Thus, $\langle v,\widetilde f\rangle$ defines a function on $\Sigma$, and
\[
\frac{\partial^2}{\partial z\,\partial\bar z}
             \langle v,\widetilde f\rangle
=g(v,D_{\frac{\partial}{\partial\bar z}}
                    \frac{\partial f}{\partial z})=0.
\]
Hence, $\langle v,\tilde{f}\rangle$ defines a constant function on $\Sigma$, and we have 
\[
g(v,\frac{\partial f}{\partial z}(p))
=g(v,\frac{\partial f}{\partial\bar z}(p))=0
\]
for all $p\in \Sigma$. In particular, we conclude that
\(
v\in V\cap V^\perp=\{0\}.
\)
It follows that $\operatorname{span}_{\mathbb R}\{a_\gamma\}=V$.
Since $V\cap\Lambda$ is discrete and spans $V$. Therefore, $V\cap\Lambda$
a lattice in $V$, and $V/(V\cap\Lambda)$ is a flat subtorus containing a translation copy of $f(\Sigma)$. Moreover, $J$ descends to this
subtorus, and $df\circ j=J\circ df$. This completes the proof.
\end{proof}

\begin{proof}[Proof of Theorem \ref{thm:flat} when $M=\RR^n$]
Take $V,J$ from Proposition \ref{prop:flat-rigidity}. For every
$x\in\Sigma$, we have
\(
f(x)-f(p)=\int_p^x df\in V.
\)
Thus $f(\Sigma)\subset f(p)+V$, and $df\circ j=J\circ df$.
\end{proof}

\section{The Partition Argument by Lee--Sidiropoulos}
In this section, we describe an energy-efficient way to partition a Riemann surface into simply connected regions. The construction is inspired by a graph-theoretic result by Lee--Sidiropoulos \cite{LeeSidiropoulos2010}. We assume $\Sigma$ is a closed, connected Riemann surface with genus $g(\Sigma)\geq 1$ and $d\mu$ a nonnegative measure. Recall that $\operatorname{sys}(\Sigma)$ is the systole of $\Sigma$; in other words, $\operatorname{sys}(\Sigma)$ is the length of the shortest non-contractible loop in $\Sigma$.

\begin{lemma}\label{lem:cutgraph}
For every $p\in\Sigma$, there exist
$4g(\Sigma)$ minimizing geodesic
segments $\alpha_1,\ldots,\alpha_{4g(\Sigma)}$ emanating
from $p$ such that $\Sigma\setminus\bigcup_{j=1}^{4g(\Sigma)}\alpha_j$ is diffeomorphic to the open unit disk in $\mathbb{C}$.
\end{lemma}

\begin{proof}
By \cite[\S 2]{Liokumovich2014},
there exist $2g(\Sigma)$ simple loops $\gamma_1,\cdots,\gamma_{2g(\Sigma)}$ with the following properties:
\begin{itemize}
    \item $p\in\gamma_j$ for each $j\in \{1,\cdots,2g(\Sigma)\}$.
    \item For every $j,k\in \{1,\cdots,2g(\Sigma)\}$ with $j\neq k$, we have $\gamma_j\cap \gamma_k=\{p\}$.
    \item For each $j\in \{1,\cdots,2g(\Sigma)\}$, $\gamma_j$ consists of two minimizing geodesic segments emanating from $p$. 
    \item $\Sigma\setminus \bigcup_{j=1}^{2g(\Sigma)}\gamma_j$ is diffeomorphic to the open unit disk in $\mathbb{C}$.
\end{itemize}
From these, the assertion follows.
\end{proof}

\begin{lemma}\label{lem:holonomy}
Let $\alpha$ be a minimizing geodesic segment in $\Sigma$, and let
\(
U=\{x\in\Sigma:d_\Sigma(x,\alpha)<\operatorname{sys}(\Sigma)/4\}.
\)
Then the inclusion $U\hookrightarrow\Sigma$
induces the trivial map on fundamental groups.
In particular, every unitary flat bundle on
$\Sigma$ admits a parallel unitary frame on $U$.
\end{lemma}

\begin{proof}
Let $\widetilde\alpha$ be a lift of $\alpha$
to the universal cover $\widetilde\Sigma$.
For every nontrivial deck transformation
$\gamma\in \pi_1(\Sigma)$ and every $x,y\in\widetilde\alpha$, the projection of
a minimizing path joining $x$ and $\gamma y$ onto $\Sigma$,
together with the corresponding subsegment of
$\alpha$, forms a noncontractible loop.
Since $\alpha$ is minimizing, the corresponding subsegment has length
at most $d_{\widetilde\Sigma}(x,\gamma y)$. The loop therefore has length
at most $2d_{\widetilde\Sigma}(x,\gamma y)$. Taking the infimum over $x,y$
gives $2d_{\widetilde\Sigma}(\widetilde\alpha,\gamma\widetilde\alpha)
\geq\operatorname{sys}(\Sigma)$.

Thus the $\frac{\operatorname{sys}(\Sigma)}{4}$-neighborhoods of distinct
lifts of $\alpha$ are disjoint, and each is mapped diffeomorphically onto $U$.
Consequently, every loop in $U$ is
contractible in $\Sigma$, this completes the proof.
\end{proof}

\begin{lemma}\label{lem:localization}
For every Lipschitz function $\psi$ on $\Sigma$,
there exist nonnegative Lipschitz functions
$\eta_0,\ldots,\eta_{4g(\Sigma)}$ with the following properties:
\begin{itemize}
    \item For each $j\in \{0,\cdots, 4g(\Sigma)\}$, the support of $\eta_j$ is contained in an open set on which every unitary flat
bundle over $\Sigma$ admits a parallel frame.
\item We have
\[
\sum_{j=0}^{4g(\Sigma)}\int_\Sigma\eta_j^2\,\psi^2\,d\mu
\geq\frac12\int_\Sigma \psi^2\,d\mu
\]
\item At almost every $x\in \Sigma$, we have
\[
\sum_{j=0}^{4g(\Sigma)}|d\eta_j(x)|^2
\leq\frac{256}{\operatorname{sys}(\Sigma)^2}\,(1+\log (4g(\Sigma)))^2.
\]
\end{itemize}

\end{lemma}

\begin{proof}
Let $S_{4g(\Sigma)}$ denote the symmetric group of  degree $4g(\Sigma)$. We consider a probability space $(0,\frac{\operatorname{sys}(\Sigma)}{4})\times S_{4g(\Sigma)}$ with the uniform probability measure $\mathbb{P}$. For each $t\in (0,\frac{\operatorname{sys}(\Sigma)}{4})$
and each permutation $\pi\in S_{4g(\Sigma)}$, we define a partition $\{P_0^{(t,\pi)},\cdots,P_{4g(\Sigma)}^{(t,\pi)}\}$ of $\Sigma$ as follows: First, we let $\alpha_1,\cdots,\alpha_{4g(\Sigma)}$  be the collection of geodesic segments constructed in Lemma \ref{lem:cutgraph}. We set
\(
P_0^{(t,\pi)}=\{x:d_\Sigma(x,\bigcup_{j=1}^{4g(\Sigma)}\alpha_j)\geq t\},
\)
and, for each $1\leq j\leq 4g(\Sigma)$, we define inductively
\(
P_j^{(t,\pi)}=
\{x:d_\Sigma(x,\alpha_{\pi(j)})<t\}
\setminus \bigcup_{k=1}^{j-1}P_k^{(t,\pi)}.
\)
Moreover, we define a collection of Lipschitz functions $\{\eta^{(t,\pi)}_j\}_{j=0}^{4g(\Sigma)}$ by
\(
\eta^{(t,\pi)}_j(x)=\min\{
1,s^{-1}\,{d_\Sigma(x,\Sigma\setminus P_j^{(t,\pi)})}
\},
\) where $s^{-1}=\frac{16}{\operatorname{sys}(\Sigma)}(1+\log (4g(\Sigma)))$.
Since the sets $\{\eta^{(t,\pi)}_j>0\}_{j=0}^{4g(\Sigma)}$ are
pairwise disjoint, we know that
\(
\sum_{j=0}^{4g(\Sigma)}|d\eta_j(x)|^2\leq s^{-2}
\) at almost every $x\in \Sigma$.

We now fix any $x\in\Sigma$, and we choose indices $i_1,\cdots,i_{4g(\Sigma)}$ such that $d_{\Sigma}(x,\alpha_{i_1})\leq \cdots\leq d_{\Sigma}(x,\alpha_{i_{4g(\Sigma)}})$.
Moreover, we assume $(t,\pi)\in (0,\frac{\operatorname{sys}(\Sigma)}{4})\times S_{4g(\Sigma)}$ such that $\sum_{j=0}^{4g(\Sigma)}\eta^{(t,\pi)}_j(x)^2<1$. In particular, we know that $\eta_0^{(t,\pi)}(x)<1$ and hence $d_{\Sigma}(x,\bigcup_{j=1}^{4g(\Sigma)}\alpha_j)\leq t+s$. We may then define $j_*$ to be the smallest element $j\in \{1,\cdots,4\,g(\Sigma)\}$ for which $d_{\Sigma}(x,\alpha_{\pi(j)})\leq t+s$. We then set $k\in \{1,\cdots,4g(\Sigma)\}$ such that $\pi(j_*)=i_k$. It then follows that
\(
d_\Sigma(x,\alpha_{\pi(j)})
>t+s
\)
for all $j<j_*$. Moreover, since $\eta_{j_*}^{(t,\pi)}(x)<1$, we must have $d_{\Sigma}(x,\alpha_{i_k})\geq t-s$ and hence $|t-d_{\Sigma}(x,\alpha_{i_k})|\leq s$. Since for every $\ell<k$, we have $d_{\Sigma}(x,\alpha_{i_\ell})\leq d_{\Sigma}(x,\alpha_{i_k})\leq t+s$, the minimality of $j_*$ implies that $\pi^{-1}(i_\ell)>j_*=\pi^{-1}(i_k)$ for all $\ell<k$. 

To summarize, we have shown that for any $x\in \Sigma$, if $\sum_{j=0}^{4g(\Sigma)}\eta_j^{(t,\pi)}(x)^2<1$, then there exists some $k\in \{1,\cdots,4g(\Sigma)\}$ such that $|t-d_{\Sigma}(x,\alpha_{i_k})|\leq s$ and $\pi^{-1}(i_\ell)>\pi^{-1}(i_k)$ for all $\ell<k$. Moreover, we observe that for any $k\in \{1,\cdots,4g(\Sigma)\}$, we have $\mathbb{P}(
|t-d_{\Sigma}(x,\alpha_{i_k})|\leq s
)\leq 8\,\operatorname{sys}(\Sigma)^{-1}s$ and
\[
\mathbb{P}\bigl(
\pi^{-1}(i_k)<\pi^{-1}(i_j)
\text{ for all }j\in\{1,\cdots,k-1\}
\bigr)=\frac1k.
\]
This implies for every $x\in\Sigma$, we have $$\mathbb{P}\Big(\sum_{j=0}^{4g(\Sigma)}\eta_j^{(t,\pi)}(x)^2<1\Big)\leq \frac{8s}{\operatorname{sys}(\Sigma)}\, \sum_{k=1}^{4\,g(\Sigma)} \frac{1}{k}\leq \frac{1}{2}.$$
Hence, we have
\(
\mathbb{E}[
\sum_{j=0}^{4g(\Sigma)}
(\eta_j^{(t,\pi)}(x))^2
]\geq\frac12
\) for every $x\in\Sigma$. Multiplying by $\psi^2$ and integrating over $\Sigma$,
we conclude that
\[
\mathbb{E}\Big[
\sum_{j=0}^{4g(\Sigma)}
\int_\Sigma
(\eta_j^{(t,\pi)}(x))^2\,\psi^2(x)\,d\mu(x)
\Big]
\geq
\frac12\int_\Sigma \psi^2(x)\,d\mu(x)
\]
Thus there exists $(t,\pi)\in(0,\frac{\operatorname{sys}(\Sigma)}{4})\times S_{4g(\Sigma)}$ such that
\[
\sum_{j=0}^{4g(\Sigma)}
\int_\Sigma
(\eta_j^{(t,\pi)}(x))^2\,\psi^2(x)\,d\mu
\geq
\frac12\int_\Sigma \psi^2(x)\,d\mu.
\]
Finally, we know that $\operatorname{supp}\eta_0^{(t,\pi)}
\Subset\Sigma\setminus \bigcup_{j=1}^{4g(\Sigma)}\alpha_j$. By Lemma \ref{lem:cutgraph}, $\Sigma\setminus \bigcup_{j=1}^{4g(\Sigma)}\alpha_j$ is diffeomorphic to the unit disk in $\mathbb{C}$, and hence every unitary flat
bundle admits a parallel unitary frame
on $\operatorname{supp}\eta_0^{(t,\pi)}$. While for each $1\leq j\leq 4g(\Sigma)$, we have
\(\operatorname{supp}\eta_j^{(t,\pi)}\Subset d_\Sigma(x,\alpha_{\pi(j)}) <\frac{\operatorname{sys}(\Sigma)}{4}\), and the claim now follows from Lemma \ref{lem:holonomy}.
\end{proof}

\section{Stable Minimal Surfaces in Manifolds with Positive Isotropic Curvature}
In this section, we assume $(M,g)$ is a Riemannian manifold of dimension at least 4 with positive isotropic curvature. Moreover, let $\Sigma$ be a closed, connected Riemann surface of genus at least 1, and let $f:\Sigma\to M$ be a conformal minimal immersion. We denote $\lambda^2=g(\frac{\partial f}{\partial z},\frac{\partial f}{\partial \bar z})$. The crucial result of this section is:

\begin{proposition}\label{prop:twisted-section}
    There exists a unitary flat complex vector bundle $E\to \Sigma$ of rank $m$, and a non-zero holomorphic section $S\in H^0(N_{\mathbb{C}}\otimes E)$ such that on a simply connected open set, if $S=\sum_{\alpha=1}^ms_{\alpha}\otimes e_{\alpha}$, then $$\sum_{\alpha=1}^mR(\frac{\partial f}{\partial z},s_{\alpha},\frac{\partial f}{\partial \bar z},\bar s_{\alpha})>0$$ away from the branch points where $|S|\neq 0$. Furthermore, if $(M,g)$ is $\kappa$-PIC, then $$\sum_{\alpha=1}^mR(\frac{\partial f}{\partial z},s_{\alpha},\frac{\partial f}{\partial \bar z},\bar s_{\alpha})\geq \frac{\kappa}{4}\,\lambda^2\,|S|^2$$ at each point of $\Sigma$.
\end{proposition}

We need the following input from algebraic geometry.

\begin{proposition}\label{prop:stable-degree-zero}
    Let $E$ be a holomorphic vector bundle of positive rank and degree zero over a compact Riemann surface $\Sigma$, then there exists a stable holomorphic vector bundle $F$ over $\Sigma$ with $\deg(F)=0$, and a non-zero holomorphic map $\varphi:F\to E$.
\end{proposition}

Proposition \ref{prop:stable-degree-zero} should be standard for experts in algebraic geometry and we give a proof of it in Appendix A. Applying it, there exists a non-zero stable holomorphic vector bundle $F$ over $\Sigma$ with $\deg(F)=0$, and a non-zero holomorphic map $\varphi:F\to N_{\mathbb C}$.
\begin{lemma}\label{lem:F-dichotomy}
One of the following holds for $F$:
\begin{enumerate}[\rm (i)]
    \item We have $g(\varphi(v),\varphi(w))=0$ for all $p\in\Sigma$ and all $v,w\in F_p$.
    \item $F$ is a holomorphic subbundle of $N_{\mathbb C}$ and $N_{\mathbb C}=F\oplus F^{\perp}$ orthogonally with respect to $g$. In particular, $F^{\perp}$ is a holomorphic subbundle of $N_{\mathbb{C}}$ of degree zero.
\end{enumerate}
\end{lemma}
\begin{proof}
    Note that the pullback of $g$ via $\varphi$ defines a holomorphic symmetric bilinear form on $F$, and hence induces a holomorphic map $\psi:F\to F^*$ defined by $\psi(v)(w)=g(\varphi(v),\varphi(w))$. If $\psi=0$, then we conclude that $g(\varphi(v),\varphi(w))=0$ for all $p\in\Sigma$ and all $v,w\in F_p$. We may then assume that $\psi\neq 0$. On the other hand, since $F$ is stable of degree zero, we know that $F^*$ is also stable of degree zero. In particular, every non-zero holomorphic map from $F$ to $F^*$ is an isomorphism. Thus, we have $F\simeq F^*$. Next, we note that the dual of $\varphi$ defines a holomorphic map $\varphi^*:N_{\mathbb C}^*\to F^*$. Since $N_{\mathbb{C}}\simeq N_{\mathbb{C}}^*$ and $F\simeq F^*$, we know that $\varphi^*$ defines a holomorphic map from $N_{\mathbb C}$ to $F$. From this, we conclude that $N_{\mathbb C}$ splits orthogonally into the direct sum of $F$ and $F^{\perp}$ with respect to $g$. This completes the proof.
\end{proof}
\begin{proposition}\label{prop:dichotomy}
    One of the following holds for $N_{\mathbb C}$:
    \begin{enumerate}[\rm (i)]
        \item There exists a non-trivial stable holomorphic vector bundle $F$ over $\Sigma$ with $\deg(F)=0$, and a non-zero holomorphic map $\varphi:F\to N_{\mathbb C}$ such that $g(\varphi(v),\varphi(w))=0$ for all $p\in\Sigma$ and all $v,w\in F_p$.
        \item $N_{\mathbb C}$ is polystable.
    \end{enumerate}
\end{proposition}

\begin{proof}
We assume that {\rm (i)} fails. By Proposition
\ref{prop:stable-degree-zero} and Lemma
\ref{lem:F-dichotomy}, we have
\(
N_{\mathbb C}=F\oplus F^\perp,
\)
where $F$ is stable of degree zero and $F^\perp$
is a holomorphic vector bundle of degree zero.

If $F^\perp=0$, then $N_{\mathbb C}=F$ is stable.
Otherwise, since the restriction of $g$ to $F^\perp$
is nondegenerate, we may apply Proposition
\ref{prop:stable-degree-zero} and Lemma
\ref{lem:F-dichotomy} to $F^\perp$. Since {\rm (i)}
fails for $N_{\mathbb C}$, it also fails for
$F^\perp$. Thus, we obtain another stable
holomorphic bundle of degree zero which splits
off orthogonally from $F^\perp$.

Repeating this construction, the rank of the remaining
orthogonal complement strictly decreases until
it is zero. Consequently, $N_{\mathbb C}$ is an
orthogonal direct sum of stable holomorphic bundles
of degree zero, and hence is polystable. This completes the proof.
\end{proof}
\begin{proposition}\label{prop:index-computation}
    Assume that $\Sigma$ is closed. Let $F$ be a unitary flat holomorphic vector bundle of rank $m$ over $\Sigma$, and let $\varphi:F\to N_{\mathbb{C}}$ be a non-zero holomorphic morphism. Then there exists a non-zero holomorphic section $S$ of $N_{\mathbb{C}}\otimes F^*$. Moreover, assume $v_1,\cdots,v_m$ is a parallel unitary frame of $F$, and let $v^1,\cdots,v^m$ be the dual frame of $v_{\alpha}$, then $S=\sum_{\alpha=1}^m\varphi(v_{\alpha})\otimes v^{\alpha}$, and
    \[
    Q_{F^*}(S,\bar S)=-4\sum_{\alpha=1}^m\int_{\Sigma}\Big(\big|D^{\tan}_{\frac{\partial}{\partial z}} (\varphi(v_{\alpha}))\big|^2+R\Big(\frac{\partial f}{\partial z},\varphi(v_{\alpha}),
        \frac{\partial f}{\partial\bar z},\overline{ \varphi(v_{\alpha})}\Big)
\Big)\,dx\,dy 
    \]
\end{proposition}
\begin{proof}
    Under $\operatorname{Hom}(F,N_{\mathbb{C}})\simeq N_{\mathbb{C}}\otimes F^*$, the morphism $\varphi$ corresponds to
$S=\sum_\alpha\varphi(v_\alpha)\otimes v^\alpha$. Since $\varphi$ is
nonzero and holomorphic, we know that $S\neq 0$ and 
\[
D^\perp_{\frac{\partial}{\partial\bar z}}S
=\sum_{\alpha=1}^m
 D^\perp_{\frac{\partial}{\partial\bar z}}\varphi(v_\alpha)\otimes v^\alpha
=0.
\]
The claim now follows from the definition of $Q_{F^*}$.
\end{proof}

\begin{proof}[Proof of Proposition \ref{prop:twisted-section}]
    Suppose first that $N_{\mathbb C}$ is polystable of degree zero. By the Narasimhan--Seshadri theorem \cite{NS}, we may choose a flat Hermitian metric $h_0$ on $N_{\mathbb C}$ compatible with its holomorphic structure. Moreover, we define a section $A$ of $\operatorname{Hom}(\overline{N}_{\mathbb{C}}, N_{\mathbb{C}})$ such that $h_0(v,Aw)=g(v,w)$ for every $v$ and $w$. Since $g$ is a holomorphic symmetric bilinear form, it is parallel for the Chern connection of $h_0$. Hence, $A$ is also parallel for the Chern connection of $h$. We then define $P$ to be the self-adjoint map from $N_{\mathbb{C}}$ to itself with $P^2=A\circ \overline{A}$. In particular, $P$ is also parallel. We then define $h(v,w)=h_0(v,Pw)$ and $\tau=P^{-1}A$, so that $h$ is a flat Hermitian metric on $N_{\mathbb{C}}$, and $\tau$ is a parallel conjugate-linear involution $\tau:N_{\mathbb C}\to N_{\mathbb C}$ with
\(
g(v,w)=h(v,\tau w)
\) for all $v$ and $w$. Therefore, at each point $p\in\Sigma$, there exists a local $h$-orthonormal frame $\{v_1,\cdots,v_{n-2}\}$ such that $g(v_{\alpha},v_{\beta})=h(v_{\alpha},v_{\beta})=\delta_{\alpha\beta}$.

We equip $E=N_{\mathbb C}^*$ with the dual unitarily flat metric, and it follows from Proposition \ref{prop:index-computation} that the identity map defines a holomorphic section $S\in H^0(N_{\mathbb C}\otimes E)$, with $S=\sum_{\alpha=1}^{n-2}v_{\alpha}\otimes v^{\alpha}$, and
\[
Q_E(S,\bar S)
\leq
-4\int_\Sigma
\sum_{\alpha=1}^{n-2}
R\Big(\frac{\partial f}{\partial z}
,v_{\alpha},\frac{\partial f}{\partial \bar z},\bar{v}_{\alpha}\Big)\,dx\,dy.
\]
At each point $p\in\Sigma$, we consider the Hermitian bilinear form over $(N_{\mathbb C})_p$ by $$H_p(v,\bar v)=R\Big(\frac{\partial f}{\partial z}
,v,\frac{\partial f}{\partial \bar z},\bar{v}\Big).$$
Since $(M,g)$ admits positive isotropic curvature, we know that away from the branch points, $H_p(v,\bar v)>0$ for all non-zero $v\in (N_{\mathbb C})_p$ with $g(v,v)=0$. By Lemma \ref{lem:average}, we conclude that $\sum_{\alpha=1}^{n-2}H_p(v_{\alpha},\bar v_{\alpha})>0$ at each $p\in\Sigma$. If $(M,g)$ is moreover $\kappa$-PIC, then the Hermitian form $H_p(v,\bar v)- \frac{\kappa}{4}\,\lambda^2\,|v|^2$ is non-negative for isotropic vector $v$. Applying Lemma \ref{lem:average} gives 
$$\sum_{\alpha=1}^{n-2}H_p(v_{\alpha},\bar v_{\alpha})\geq \frac{\kappa}{4}\lambda^2\sum_{\alpha=1}^{n-2}|v_{\alpha}|^2=\frac{\kappa}{4}\lambda^2\,|S|^2.$$ 

Suppose now that $N_{\mathbb{C}}$ is not polystable, then by Proposition \ref{prop:dichotomy}, there exists a non-trivial stable holomorphic vector bundle $F$ over $\Sigma$ with $\deg(F)=0$, and a non-zero holomorphic map $\varphi:F\to N_{\mathbb C}$ such that $g(\varphi(v),\varphi(w))=0$ for all $p\in\Sigma$ and all $v,w\in F_p$.
By the Narasimhan--Seshadri theorem, $F$ admits a
unitary flat connection compatible with its holomorphic
structure. The map $\varphi$ then defines a non-zero
holomorphic section
\(
S\in H^0(N_{\mathbb C}\otimes F^*),
\) such that if $v_1,\cdots,v_m$ is a local parallel unitary frame of $F$, then $S=\sum_{\alpha=1}^m\varphi(v_{\alpha})\otimes v^{\alpha}$. Moreover, since $(M,g)$ has
positive isotropic curvature, we have
\[
R\Big(\frac{\partial f}{\partial z},\varphi(v),\frac{\partial f}{\partial \bar z},
\overline{\varphi(v)}\Big)>0
\]
away from the branch points whenever $\varphi(v)\neq0$. If $(M,g)$ is moreover $\kappa$-PIC, then we have 
\[
R\Big(\frac{\partial f}{\partial z},\varphi(v),\frac{\partial f}{\partial \bar z},
\overline{\varphi(v)}\Big)\geq\frac{\kappa}{4}\lambda^2\,|\varphi(v)|^2.
\]
Hence, we conclude that
\[
\sum_{\alpha=1}^m
R\Big(
\frac{\partial f}{\partial z},\varphi(v_\alpha),
\frac{\partial f}{\partial\bar z},
\overline{\varphi(v_\alpha)}
\Big)
\geq
\frac{\kappa}{4}\lambda^2
\sum_{\alpha=1}^m|\varphi(v_\alpha)|^2
=
\frac{\kappa}{4}\lambda^2|S|^2.
\]
This completes the proof.
\end{proof}

\begin{proof}[Proof of Theorem \ref{thm:pic}]
If $\Sigma$ has genus zero, the immersion is unstable by
Micallef--Moore \cite{MM1988}. Otherwise if $g(\Sigma)\geq 1$. Proposition \ref{prop:twisted-section} implies there exists a unitary flat complex vector bundle $E\to \Sigma$ , and a non-zero holomorphic section $S\in H^0(N_{\mathbb{C}}\otimes E)$ such that $Q_E(S,\bar{S})<0$. Hence, by Corollary \ref{cor:destablization}, $f$ cannot be covering stable. This completes the proof.
\end{proof}
\begin{proof}[Proof of Theorem \ref{thm:PICsystole}]
    We assume that $(M,g)$ is $\kappa$-PIC for some $\kappa>0$. By Proposition \ref{prop:twisted-section}, there exists a unitary flat complex vector bundle $E\to \Sigma$ of rank $m$, and a non-zero holomorphic section $S\in H^0(N_{\mathbb{C}}\otimes E)$ such that on a simply connected open set, if $S=\sum_{\alpha=1}^ms_{\alpha}\otimes e_{\alpha}$, then 
    $$\sum_{\alpha=1}^mR(\frac{\partial f}{\partial z},s_{\alpha},\frac{\partial f}{\partial \bar z},\bar s_{\alpha})\geq \frac{\kappa}{4}\,\lambda^2\,|S|^2$$ at each point of $\Sigma$. Next, we apply Lemma \ref{lem:localization} with $\psi=|S|$, it then follows that there exists a family of non-negative Lipschitz functions $\eta_0,\cdots,\eta_{4g(\Sigma)}$ such that for each $j$, we can choose a parallel unitary frame $e^{(j)}_1,\cdots,e^{(j)}_m$ on an open neighborhood of $\operatorname{supp}\eta_j$. Moreover, we may write $S=\sum_{\alpha=1}^ms^{(j)}_{\alpha}\otimes e_{\alpha}^{(j)}$ on $\operatorname{supp}\eta_j$, where each $s^{(j)}_{\alpha}$ is holomorphic. In particular, each $\eta_j\,s^{(j)}_{\alpha}$ defines a Lipschitz section of $N_{\mathbb{C}}$. Since $f:\Sigma \to M$ is stable, and by a standard approximation argument, we have
    \begin{align*}
        0&\leq \sum_{\alpha=1}^mQ(\eta_j\,s^{(j)}_{\alpha},\eta_j\,\bar s^{(j)}_{\alpha})\\
        &\leq 4\sum_{\alpha=1}^m\int_{\Sigma}\Big(|D^{\perp}_{\frac{\partial}{\partial \bar z}}(\eta_j\,s^{(j)}_{\alpha})|^2-\frac{\kappa}{4}\lambda^2\eta_j^2\,|s_{\alpha}^{(j)}|^2\Big)\,dx\,dy\\
        &\leq 4\int_{\Sigma}\Big|\frac{\partial \eta_j}{\partial\bar z}\Big|^2\,|S|^2\,dx\,dy-\kappa\int_{\Sigma}\lambda^2\,\eta_j^2\,|S|^2\,dx\,dy\\
        &=\int_{\Sigma}\Big(|d\eta_j|^2-\frac{\kappa}{2}\eta_j^2\Big)\,|S|^2\,d\mu.
    \end{align*}
    Here, $d\mu=2\lambda^2\,dx\,dy$ is the volume form of $\Sigma$ with respect to the metric $f^*g$. Summing over $j$ and using the estimates in Lemma \ref{lem:localization} gives
\begin{align*}
\frac{\kappa}{4}\int_\Sigma |S|^2\,d\mu
&\leq
\frac{\kappa}{2}
\sum_{j=0}^{4g(\Sigma)}\int_\Sigma\eta_j^2|S|^2\,d\mu\\
&\leq
\sum_{j=0}^{4g(\Sigma)}\int_\Sigma|d\eta_j|^2|S|^2\,d\mu\\
&\leq
\frac{256}{\operatorname{sys}(\Sigma)^2}
\bigl(1+\log(4g(\Sigma))\bigr)^2
\int_\Sigma |S|^2\,d\mu.
\end{align*}
Since $\int_{\Sigma}|S|^2\,d\mu>0$, the theorem follows directly.
\end{proof}

\section{Stable Minimal Immersions with Logarithmic Systole Growth}\label{sec:log-systole}
In this section, we construct an example of a stable minimal surface of genus two in $M=(S^3\times S^1)\#(S^3\times S^1)$ with a positive isotropic curvature metric $g$ on $M$, such that there exists a sequence of stable covers whose systole grows logarithmically in the genus. In particular, this shows that the logarithmic dependence in Theorem \ref{thm:PICsystole} is sharp.
\begin{theorem}\label{thm:logarithmic-systole}
There exist a constant $\kappa>0$, a $\kappa$-PIC metric $g$ on $M$, a totally geodesic embedding $f:\Sigma\to M$ of a closed oriented surface of genus two, and a sequence of connected finite regular covers $\pi_k:\Sigma_k\to\Sigma$ with the following properties.
\begin{enumerate}[\rm (i)]
\item For $f$ and for each $f\circ\pi_k$ the index form satisfies $Q(s,\bar s)\geq\frac{\kappa}{2}\int |s|^2\,d\mu$ for every smooth complex normal section $s$.
\item $g(\Sigma_k)\to\infty$ and there exists some constant $a>0$ such that 
$\operatorname{sys}(\Sigma_k)\geq a\log g(\Sigma_k)$ for each $k\geq 1$. Here, the systole of $\Sigma_k$ is computed with respect to $(f\circ\pi_k)^*g$.
\end{enumerate}
\end{theorem}

The proof of Theorem \ref{thm:logarithmic-systole} is divided into three steps. We first construct an $O(2)$-invariant locally conformally flat PIC metric $g$ on $(S^3\times S^1)\#(S^3\times S^1)$ whose fixed points set $\Sigma$ is a minimal surfaces with a flat normal bundle with nontrivial holonomy. In step 2, we construct conformal deformations on $M$, keeping the induced metric on $\Sigma$ and its normal connection, while decreasing the ambient PIC term in the normal direction of $\Sigma$. In step 3, we show that suitable high-degree covers of $\Sigma$ is stable: high-degree random covers any genus two surface has a uniformly positive first Laplacian eigenvalue that dominates the small positive curvature term appeared in step 2.

\begin{proposition}\label{prop:example-background}
There exists a locally conformally flat metric $\bar{g}$ of positive scalar curvature on $M$ that is invariant under an $O(2)$-action whose common fixed point set is an embedded totally geodesic surface $\Sigma$ of genus two. Moreover, the real normal bundle splits as $N \Sigma=L\oplus L$ for some nontrivial flat real line bundle $L$ with its induced connection.
\end{proposition}
\begin{proof}
    Let $S^4$ be the round unit sphere in $\mathbb R^3\times\mathbb R^2$. For $A\in\mathrm O(2)$, set $\Psi_A(x,y)=(x,Ay)$ for $(x,y)\in S^4$. The fixed point set is the round sphere $S^2=S^4\cap(\mathbb R^3\times\{0\})$. Fix four distinct points $p_1^\pm,p_2^\pm$ on $S^2$. By the Gromov--Lawson surgery construction, there exists an $O(2)$-invariant metric $\bar{g}$ conformal to the round metric on $S^4 \setminus \{p_1^+,p_1^-,p_2^+,p_2^-\}$ with uniformly positive scalar curvature and a round cylindrical end of small radius $\varepsilon$ at each puncture.

    For $j=1,2$, we glue the cylindrical ends near $p_j^+$ and $p_j^-$ as follows. Consider the orientation-reversing linear maps $F_1,F_2$ on $\mathbb R^4$:
    \begin{align*}
        &F_1(x_1,x_2,x_3,x_4)=(-x_1,x_2,x_3,x_4),\\
        &F_2(x_1,x_2,x_3,x_4)=(-x_1,x_2,-x_3,-x_4).
    \end{align*}
    On each cylindrical end, choose coordinates $(\omega,r)$ compatible with the orientation of $S^4$. In these coordinates, we have $\bar{g}=dr\otimes dr+\varepsilon^2g_{S^3}$, where $g_{S^3}$ is the unit round metric on $S^3$. At each puncture $p\in \{p_1^+,p_1^-,p_2^+,p_2^-\}$, with $\omega$ parametrizes the unit sphere in $T_pS^4$, and for $A\in O(2)$, we have $\Psi_A(\omega,r)=(\omega_1,\omega_2,A(\omega_3,\omega_4),r)$. Next, for $j\in \{1,2\}$, we identify $(\omega,0)$ near $p_j^+$ with $(F_j(\omega),0)$ near $p_j^-$. Since $F_j$ is orientation reversing, the resulting manifold is $(S^3\times S^1) \# (S^3\times S^1) $. Since the map $(\omega,r)\to (F_j(\omega),-r)$ is an isometry of the round cylinder, we conclude that $\bar{g}$ is smooth and locally conformally flat with positive scalar curvature. Moreover, $\Psi_A$ commutes with both $F_1$ and $F_2$. Hence, the common fixed point set on each end is parametrized by $\omega_3=\omega_4=0$. The fixed surface $\Sigma$ is therefore obtained by removing four open disks from $S^2$ and identifying their boundary circles into two pairs. We conclude that $\Sigma$ is totally geodesic with genus two.

    On $S^2$ with the four disks removed, the normal bundle is spanned by $(0,0,0,1,0)$ and $(0,0,0,0,1)$. The differential of $\Psi_{\operatorname{diag}(1,-1)}$ fixes the first vector and reverses the second vector. Since the normal parallel transport commutes with $\Psi_{\operatorname{diag}(1,-1)}$, so it preserves both normal lines. For $j=1,2$, let $\gamma_j$ be a loop that crosses the $j$-th gluing circle once and avoids the other. The first gluing leaves $(\omega_3,\omega_4)$ unchanged, and the second sends $(\omega_3,\omega_4)$ to $(-\omega_3,-\omega_4)$. Hence, we have $\operatorname{Hol}_{\gamma_1}=\id$ and $\operatorname{Hol}_{\gamma_2}=-\id$. Consequently, we have $N\Sigma=L\oplus L$ for some nontrivial flat normal bundle $L$, with holonomy $-1$ along $\gamma_2$.
\end{proof}

\begin{proposition}\label{prop:example-scalar}
Let $\bar{g}$ and $\Sigma\subset M$ be as in Proposition \ref{prop:example-background}. For every $0<3\kappa<\min_M\mathrm{scal}_{\bar g}$, there exists an
$\mathrm O(2)$-invariant metric $g$ on $M$ conformal to $\bar g,$ such that for all $x\in \Sigma$, we have $\bar{g}(x)=g(x)$ and $\bar{D}g(x)=0$. Here, $ \bar D$ is the Levi--Civita connection of $\bar g$. Moreover, we have
$\mathrm{scal}_g\geq 3\kappa$ on $M$ and
$\mathrm{scal}_g=3\kappa$ at each point of $\Sigma$.
\end{proposition}

\begin{proof}
We choose some $r_0>0$ sufficiently small, and an $O(2)$ invariant smooth function $\rho:M\to[0,\infty)$, with $\rho^{-1}(0)=\Sigma$
and $\rho(x)=d_{\bar g}(x,\Sigma)^2$ whenever
$d_{\bar g}(x,\Sigma)<r_0$.
We have $d\rho=0$ and $\Delta_{\bar g}\rho=4$ at each point of $\Sigma$.
By Proposition \ref{prop:example-background}, $\Psi_{-I}$ fixes
$\Sigma$ and reverses the normal vectors. Thus the two invariant
functions $\Delta_{\bar g}\rho-4$ and
$\langle d\,\mathrm{scal}_{\bar g},d\rho\rangle_{\bar g}$
are even in normal exponential coordinates and vanish along $\Sigma$.
Hence, there exists some $C>0$ such that
\[
|\Delta_{\bar g}\rho-4|
+|\langle d\,\mathrm{scal}_{\bar g},d\rho\rangle_{\bar g}|
\leq C\rho
\]
at each point of $M$. For $0<\varepsilon\leq1$, we set $g=w^2\bar g$ with
\[
w=1+\frac{\varepsilon}{24}
(1-e^{-\frac{\rho}{\varepsilon}})(\mathrm{scal}_{\bar g}-3\kappa).
\]
Since
$\rho e^{-\frac{\rho}{\varepsilon}}\leq
\varepsilon(1-e^{-\frac{\rho}{\varepsilon}})$, we have
\begin{align*}
6\Delta_{\bar g}w
&=\frac{1}{4}e^{-\frac{\rho}{\varepsilon}}
(\mathrm{scal}_{\bar g}-3\kappa)
\Big(\Delta_{\bar g}\rho-\frac{|d\rho|_{\bar g}^2}{\varepsilon}\Big)\\
&\quad+\frac{1}{2}e^{-\frac{\rho}{\varepsilon}}
\langle d\,\mathrm{scal}_{\bar g},d\rho\rangle_{\bar g}
+\frac{\varepsilon}{4}(1-e^{-\frac{\rho}{\varepsilon}})
\Delta_{\bar g}\mathrm{scal}_{\bar g}\\
&\leq e^{-\frac{\rho}{\varepsilon}}(\mathrm{scal}_{\bar g}-3\kappa)
+C\varepsilon(1-e^{-\frac{\rho}{\varepsilon}}).
\end{align*}
Since $1\leq w\leq1+\frac{1}{24}\max_M(\mathrm{scal}_{\bar g}-3\kappa)$, we also have $w^3-w=w(w+1)(w-1)\leq  C\varepsilon(1-e^{-\frac{\rho}{\varepsilon}})$. Consequently, we have
\begin{align*}
w^3(\mathrm{scal}_g-3\kappa)
&=w(\mathrm{scal}_{\bar g}-3\kappa)-3\kappa(w^3-w)-6\Delta_{\bar g}w\\
&\geq(1-e^{-\frac{\rho}{\varepsilon}})
(\mathrm{scal}_{\bar g}-3\kappa-C\varepsilon)
\geq 0,
\end{align*}
provided $\varepsilon$ is sufficiently small so that $C\varepsilon<\min_M\mathrm{scal}_{\bar g}-3\kappa$.

At each point of $\Sigma$, we have $w=1$ and $dw=0$. Hence we have $g(x)=\bar g(x)$ and $\bar Dg=2w\,dw\otimes\bar g$ implies
$\bar Dg(x)=0$ for all $x\in \Sigma$.
The formula for $6\Delta_{\bar g}w$ gives
$6\Delta_{\bar g}w=\mathrm{scal}_{\bar g}-3\kappa$ on $\Sigma$,
so $\mathrm{scal}_g=\mathrm{scal}_{\bar g}-6\Delta_{\bar g}w=3\kappa$
at each point of $\Sigma$. This completes the proof.
\end{proof}
\begin{lemma}\label{lem:example-covers}
    Let $(\Sigma,h)$ be a closed oriented Riemannian surface of genus at
least two, and let $L\to\Sigma$ be a
nontrivial flat real line bundle with metric connection $D^L$. There exist constants $a,b>0$ and a sequence of
connected finite regular covers $\pi_k:\Sigma_k\to\Sigma$ of degree tending to infinity with the following properties.
\begin{itemize}
\item For every $k\geq1$, we have
\(
\operatorname{sys}_{\pi_k^*h}(\Sigma_k)\geq a\log g(\Sigma_k).
\) Here, the systole is computed with respect to the metric $h_k=\pi_k^*h$.
\item For every $k\geq 0$, and every smooth section $s$ of $\pi_k^*L\otimes_{\mathbb{R}} \mathbb{C}$, we have
\[
\int_{\Sigma_k}|D^Ls|^2_{h_k}\,d\mu_{h_k}\geq b\int_{\Sigma_k}|s|^2\,d\mu_{h_k}
\]
\end{itemize}
\end{lemma}
\begin{proof}
By a result of Breuillard and Gamburd \cite[Corollary 1.2]{BreuillardGamburd2010}, there exist $c>0$, and a sequence of odd primes $p_k\to\infty$ such that every Cayley graph of $\mathrm{SL}_2(\mathbb F_{p_k})$ is a $c$-expander: equivalently, the first non-zero eigenvalue of the graph Laplacian is uniformly bounded away from zero. We now fix a hyperbolic metric $g_{\operatorname{hyp}}$ with $h\geq c_0^2\,g_{\operatorname{hyp}}$ for some $c_0>0$. By the recent work of Liu and Petri \cite[Theorem 3.1 and Section 3.6]{LiuPetri2026}, after passing to a subsequence of $p_k$, for each $k$ there exists a surjective homomorphism $\rho_k:\pi_1(\Sigma)\to \mathrm{SL}_2(\mathbb F_{p_k})$ such that the associated cover $\pi_k:\Sigma_k\to \Sigma$ satisfies
\(
\operatorname{sys}_{\pi_k^*h}(\Sigma_k)
\geq c_0\,\operatorname{sys}_{\pi_k^*g_{\operatorname{hyp}}}(\Sigma_k)
\geq a\log g(\Sigma_k).
\) Here, $a>0$ is independent of $k$. Moreover, we note that $\deg(\pi_k)=|\mathrm{SL}_2(\mathbb F_{p_k})|=p_k(p_k^2-1)$ which tends to infinity as $k\to\infty$.

Next, we set $\hat{\Sigma}_0=\{(x,v):x\in\Sigma,v\in L_x,|v|=1\}$ and set $\hat{\pi}_0(x,v)=x$. Since $L$ is non-trivial, $\hat{\pi}_0$ is a connected double cover. For each $k\geq 1$, we write $\hat{\pi}_k:\hat{\Sigma}_k\to \Sigma_k$ be its pullback by $\pi_k$. Next, we observe that $\rho_k\circ(\widehat\pi_0)_*:\pi_1(\hat{\Sigma}_0)\to \mathrm{SL}_2(\mathbb F_{p_k})$ for each $k\geq 1$. Indeed, the group $\mathrm{SL}_2(\mathbb F_{p_k})$ is generated by unipotent elements of odd order, so it does not contain quotient of order two. Since the image of $\rho_k\circ(\widehat\pi_0)_*$ has index at most two, we conclude that it must be the whole group. Hence, $\widehat\Sigma_k$ is a connected regular $SL_2(\mathbb{F}_p)$-cover of $\hat{\Sigma}_0$.

We then write $\hat{h}_k=\hat{\pi}_k^* h_k$. By the classical Brooks's theorem \cite[Theorem 1]{Brooks1986} applying to the Cayley graphs defined by the images of a fixed finite symmetric generating set of $\pi_1(\widehat\Sigma_0)$ under $\rho_k\circ(\widehat\pi_0)_*$, we conclude that there exists some $b>0$ independent of $k$, so that  $\lambda_1(\widehat\Sigma_k,\widehat h_k)\geq b>0$ for all $k$. 

For $(x,v)\in \hat{\Sigma}_k$, where $v$ is a unit vector in $L_{\pi_k(x)}$, we write $s(x)=u(x,v)\,v$ and it follows that $u(x,v)=-u(x,-v)$. Since the map $(x,v)\to (x,-v)$ preserves $\hat{h}_k$, we conclude that $\int_{\hat{\Sigma}_k}u\,d\mu_{\hat{h}_k}=0$. From this, we conclude that
\[
\int_{\Sigma_k}|D^Ls|_{h_k}^2\,d\mu_{h_k}=\frac{1}{2}\int_{\hat{\Sigma}_k}|du|_{\hat{h}_k}^2\,d\mu_{\hat{h}_k}\geq \frac{b}{2}\int_{\hat{\Sigma}_k}|u|^2\,d\mu_{\hat{h}_k}=b\int_{\Sigma_k}|s|^2\,d\mu_{h_k}.
\]
This completes the proof.
\end{proof}
\begin{proof}[Proof of Theorem \ref{thm:logarithmic-systole}]
Let $\bar g$, $f:\Sigma\to M$, and
$L\to\Sigma$ be as in Proposition
\ref{prop:example-background}. Put $h=f^*\bar g$ and apply
Lemma \ref{lem:example-covers} to the line bundle $L$. Choose $0<\kappa\leq b$ with
$3\kappa<\min_M\mathrm{scal}_{\bar g}$, and take $g$ from Proposition \ref{prop:example-scalar}. Then $f$ is totally geodesic with $f^*g=h$, and the normal bundle is still given by $L\oplus L$. Moreover, we have a sequence of connected regular covers $\pi_k:\Sigma_k\to\Sigma$ such that $g(\Sigma_k)=1+\deg(\pi_k)\to\infty$ and with $\operatorname{sys}(\Sigma_k)\geq a\,\log g(\Sigma_k)$ for each $k\geq 1$.

Let $R$ and $\text{\rm Ric}$ denote the Riemann curvature tensor and the Ricci curvature of $g$.
In a $g$-orthonormal four-frame $e_1,\ldots,e_4$, write
$R_{ijk\ell}=R(e_i,e_j,e_k,e_\ell)$ and
$\text{\rm Ric}_{ij}=\text{\rm Ric}(e_i,e_j)$.
Since $g$ is locally conformally flat, we have $R_{1234}=0$ and
\[
R_{ijij}=\frac{1}{2}(\text{\rm Ric}_{ii}+\text{\rm Ric}_{jj})
-\frac{1}{6}\mathrm{scal}_g
\]
for $i\neq j$. Consequently, we have
\begin{align*}
&R_{1313}+R_{1414}+R_{2323}+R_{2424}-2R_{1234}=\frac{1}{3}\mathrm{scal}_g\geq\kappa.
\end{align*}
Thus $g$ is $\kappa$-PIC.

Fix an orthonormal $\{e_1,e_2\}$ of $\Sigma$, then the bilinear form on $N\Sigma$ defined by $(X,Y)\mapsto R(e_1,X,e_1,Y)+R(e_2,X,e_2,Y)$ is invariant under the $O(2)$-action. Evaluating the trace implies that $R(e_1,X,e_1,Y)+R(e_2,X,e_2,Y)=\frac{\kappa}{2}g(X,Y)$. 

We fix $k\geq 0$ and let $s$ be a smooth complex normal vector field along $f\circ \pi_k$. The complexified normal bundle of $f\circ \pi_k$ is given by $(\pi_k^*L\otimes_{\mathbb{R}} \mathbb{C})\oplus (\pi_k^*L\otimes_{\mathbb{R}} \mathbb{C})$. We write $D$ for the induced normal connection. Since $f\circ\pi_k$ is totally geodesic, we have
\begin{align*}
Q(s,\bar s)&=\int_{\Sigma_k}\Big(|Ds|_{h_k}^2-\frac{\kappa}{2}|s|^2\Big)
\,d\mu_{h_k}\\
&\geq \Big(b-\frac{\kappa}2{}\Big)\int_{\Sigma_k}|s|^2\,d\mu_{h_k}\geq \frac{\kappa}{2}\int_{\Sigma_k}|s|^2\,d\mu_{h_k}.
\end{align*}
This completes the proof.
\end{proof}

\appendix

\section{Facts about Holomorphic Vector Bundles over a Riemann Surface}
The goal of this appendix is to prove Proposition \ref{prop:stable-degree-zero}.  We start by reviewing some basics of holomorphic vector bundles over a Riemann surface $\Sigma$. We use $\cO_\Sigma$ to denote the sheaf of holomorphic functions. It is well-known that the category of holomorphic vector bundles over $\Sigma$ is equivalent to the category of finite rank locally free sheaves as $\cO_\Sigma$-modules \cite[Ch.~II, \S3.1, pp.~84--86]{Demailly2012}. Under this isomorphism, holomorphic bundle maps correspond to $\cO_\Sigma$-linear morphisms of sheaves. Consequently, we use the same notation for a holomorphic vector bundle as its sheaves of holomorphic sections.

We start with the following Lemma.
\begin{lemma}\label{lem:lowering-degree}
Let $G$ be a holomorphic vector bundle over a compact Riemann surface $\Sigma$ with $\deg(G)\geq 0$, then there exists a locally free subsheaf $\hat{G}\subset G$ over $\Sigma$ such that $\deg(\hat{G})=0$ and $\operatorname{rank}(\hat{G})=\operatorname{rank}(G)$.    
\end{lemma}

\begin{proof}
If $\deg(G)=0$, we take $\hat{G}=G$. Otherwise, we perform elementary modifications to $G$ at points of $\Sigma$ as follows \cite[\S4.4, p.~179]{Ramanan2009}. Fix $p\in \Sigma$ and a non-zero linear function $\ell$ on the fiber $G_p$. Define a new sheaf $G'$ (as the result of an elementary modification to $G$) via the exact sequence of sheaves
\[0\longrightarrow G' \longrightarrow G\xrightarrow{\ell \circ \ev_p} \mathbb{C}_p\longrightarrow 0, \]
where $\mathbb{C}_p$ denotes the skyscraper sheaf supported at $p$. Equivalently, $G'$ consists of local holomorphic sections $s$ of $G$ satisfying $\ell(s(p))=0$ near $p$, and impose no condition away from $p$. Precisely, in a local coordinate $z$ centered around $p$, we may choose a holomorphic frame $e_1,\cdots,e_r$ of $G$ and assume that $\ell$ is the first coordinate at $p$. Then $ze_1,e_2,\cdots,e_r$ generate $G'$ freely around $p$, and $G'=G$ away from $p$. Thus $\rank (G')=\rank (G)$ and $G'$ is also locally free. Moreover, we have $\det(G')\simeq \det(G)\otimes \mathcal{O}_{\Sigma}(-p)$ and hence $\deg (G') = \deg (G) -1$. 

By performing $\deg(G)$ elementary
modifications at points of $\Sigma$, we obtain a locally
free subsheaf $\hat G\subset G$ of the same rank with
\(
\deg(\hat G)=\deg(G)-\deg(G)=0.
\) This completes the proof.
\end{proof}

We now give a proof of Proposition \ref{prop:stable-degree-zero}. We let 
$$\mathcal{M}=\{\text{non-zero locally free subsheaves of $E$ with degree zero}\}.$$
Clearly, $E\in\mathcal{M}$ and hence $\mathcal{M}$ is non-empty. We let $F\in \mathcal{M}$ such that $F$ has minimal rank among all elements of $\mathcal{M}$. We now claim that $F$ is stable. Suppose to the contrary, that there exists a proper non-zero coherent subsheaf $G\subset F$ such that $0<\operatorname{rank}(G)<\operatorname{rank}(F)$ and $\mu(G)\geq \mu(F)=0$. In particular, we have $\deg(G)\geq 0$. Moreover, since $G$ is a torsion-free coherent sheaf over a Riemann surface, we must have $G$ is locally free and hence a holomorphic vector bundle \cite[Ch.~V, \S5]{Kobayashi1987}. By Lemma \ref{lem:lowering-degree}, we can find a locally free subsheaf $\hat{G}\subset G$ over $\Sigma$ such that $\deg(\hat{G})=0$ and $\operatorname{rank}(\hat{G})=\operatorname{rank}(G)$. Therefore, $\hat{G}$ is a non-zero locally free subsheaf of $E$ with $\deg(\hat{G})=0$ and $\operatorname{rank}(\hat{G})<\operatorname{rank}(F)$. This contradicts the minimality of $F$.

Hence, $F$ is a stable degree zero locally free subsheaf of $E$. We conclude that $F$ is a holomorphic vector bundle, and the inclusion map $\varphi:F\to E$ is non-zero. This completes the proof of Proposition \ref{prop:stable-degree-zero}. 

\begin{remark}
    An example to keep in mind is when $p\in \Sigma$ and $E = \cO(p)\oplus \cO(-p)$. In this case, we perform the elementary modification on the $\cO(p)$ summand and obtain $F = \cO_\Sigma$. The holomorphic mapping $\varphi: F\to E$ is given explicitly as
    \[u \mapsto (us, 0) , \]
    where $s$ is the canonical nonzero section of $\cO(p)$ whose zero divisor is $p$.

    In general, the above construction also guarantees that $\varphi$ is fiberwise injective away from finitely many points.
\end{remark}

\section{A Statement in Complex Linear Algebra}

\begin{lemma}\label{lem:average}
Let $Z$ be a complex vector space of dimension $m\geq2$,
equipped with a nondegenerate symmetric complex bilinear
form $g$ and a Hermitian inner product $h$. Suppose
there exists a basis $v_1,\ldots,v_m$ satisfying
\(
g(v_\alpha,v_\beta)
=h(v_\alpha,v_\beta)=\delta_{\alpha\beta}.
\)
Let $H$ be a Hermitian form on $Z$ such that there exists some $\kappa\geq 0$ and
\(
H(v,\bar v)\geq \kappa\,h(v,\bar{v})
\)
for every nonzero $v\in Z$ satisfying $g(v,v)=0$.
Then
\(
\sum_{\alpha=1}^m H(v_\alpha,\bar v_\alpha)\geq m\,\kappa.
\)
\end{lemma}

\begin{proof}
For every $\alpha<\beta$, we have
\[
0=g(v_\alpha+i v_\beta,v_\alpha+i v_\beta)=g(v_\alpha-i v_\beta,v_\alpha-i v_\beta).
\]
Hence, it follows that
\begin{align*}
    4\kappa&\leq H(v_\alpha+i v_\beta,\bar v_\alpha-i \bar v_\beta)+H(v_\alpha-i v_\beta,\bar v_\alpha+i \bar v_\beta)\\
    &=2H(v_\alpha,\bar v_\alpha)+2H(v_{\beta},\bar v_{\beta}).
\end{align*}
Summing over all $\alpha<\beta$, we obtain
\(
(m-1)\sum_{\alpha=1}^m
H(v_\alpha,\bar v_\alpha)\geq m(m-1)\kappa.
\)
This proves the assertion.
\end{proof}

\bibliographystyle{amsplain}
\bibliography{bib}

\providecommand{\bysame}{\leavevmode\hbox to3em{\hrulefill}\thinspace}
\providecommand{\MR}{\relax\ifhmode\unskip\space\fi MR }
\providecommand{\MRhref}[2]{%
  \href{http://www.ams.org/mathscinet-getitem?mr=#1}{#2}
}
\providecommand{\href}[2]{#2}
\begin{thebibliography}{10}

\bibitem{AS}
L.~V. Ahlfors and L.~Sario, \emph{{Riemann Surfaces}}, Princeton Mathematical
  Series, vol.~26, Princeton University Press, Princeton, NJ, 1960.

\bibitem{ArezzoMicallef2000}
C.~Arezzo and M.~J. Micallef, \emph{{Minimal surfaces in flat tori}}, Geom.
  Funct. Anal. \textbf{10} (2000), 679--701.

\bibitem{AMP2002}
C.~Arezzo, M.~J. Micallef, and G.~P. Pirola, \emph{{Stable minimal surfaces of
  finite total curvature}}, Comm. Anal. Geom. \textbf{10} (2002), no.~1,
  11--22.

\bibitem{BreuillardGamburd2010}
Emmanuel Breuillard and Alex Gamburd, \emph{Strong uniform expansion in
  $\mathrm{SL}(2, p)$}, Geometric and Functional Analysis \textbf{20} (2010),
  no.~5, 1201--1209.

\bibitem{Brooks1986}
Robert Brooks, \emph{The spectral geometry of a tower of coverings}, Journal of
  Differential Geometry \textbf{23} (1986), no.~1, 97--107.

\bibitem{CMR24}
G.~Catino, P.~Mastrolia, and A.~Roncoroni, \emph{{Two rigidity results for
  stable minimal hypersurfaces}}, Geom. Funct. Anal. \textbf{34} (2024), no.~1,
  1--18.

\bibitem{CKM2025}
D.~R. Cheng, S.~Karigiannis, and J.~Madnick, \emph{{Holomorphicity of parabolic
  stable minimal surfaces of high codimension}}, arXiv:2509.22155, 2025,
  revised 2026.

\bibitem{CL23}
O.~Chodosh and C.~Li, \emph{{Stable anisotropic minimal hypersurfaces in
  $\mathbf{R}^4$}}, Forum Math. Pi \textbf{11} (2023), Paper No.~e3, 22~pp.

\bibitem{CL24}
\bysame, \emph{{Stable minimal hypersurfaces in $\mathbf{R}^4$}}, Acta Math.
  \textbf{233} (2024), no.~1, 1--31.

\bibitem{CLMS26}
O.~Chodosh, C.~Li, P.~Minter, and D.~Stryker, \emph{{Stable minimal
  hypersurfaces in $\mathbf{R}^5$}}, Ann. of Math. (2) \textbf{204} (2026),
  no.~1, 423--453.

\bibitem{Demailly2012}
J.-P. Demailly, \emph{{Complex Analytic and Differential Geometry}}, Institut
  Fourier, Universit\'e de Grenoble~I, online manuscript, 2012, version of
  June~21, 2012.

\bibitem{dCP79}
M.~do~Carmo and C.~K. Peng, \emph{{Stable complete minimal surfaces in
  $\mathbf{R}^3$ are planes}}, Bull. Amer. Math. Soc. (N.S.) \textbf{1} (1979),
  no.~6, 903--906.

\bibitem{EFM26}
Nick Edelen, Luis~Atzin Franco~Reyna, and Paul Minter, \emph{Entire
  area-minimizing surfaces in $\mathbf{R}^n$ of density two are quadratic or
  planar}, \url{https://arxiv.org/abs/2606.15795} (2026).

\bibitem{FCS80}
D.~Fischer-Colbrie and R.~Schoen, \emph{{The structure of complete stable
  minimal surfaces in $3$-manifolds of nonnegative scalar curvature}}, Comm.
  Pure Appl. Math. \textbf{33} (1980), no.~2, 199--211.

\bibitem{Fraser2026}
A.~Fraser, \emph{{Stability theory for minimal surfaces}}, {Proceedings of the
  International Congress of Mathematicians 2026}, vol.~4, SIAM, 2026,
  pp.~139--159.

\bibitem{FraserSchoen2025}
A.~Fraser and R.~Schoen, \emph{{Stability and largeness properties of minimal
  surfaces in higher codimension}}, J. Eur. Math. Soc. (2025).

\bibitem{Fraser2003}
A.~M. Fraser, \emph{{Fundamental groups of manifolds with positive isotropic
  curvature}}, Ann. of Math. (2) \textbf{158} (2003), no.~1, 345--354.

\bibitem{HS85}
R.~Hardt and L.~Simon, \emph{{Area minimizing hypersurfaces with isolated
  singularities}}, J. Reine Angew. Math. \textbf{362} (1985), 102--129.

\bibitem{HLW26}
H.~Hong, H.~Li, and G.~Wang, \emph{{The stable Bernstein theorem in
  $\mathbb{R}^7$}}, preprint, arXiv:2609.15720, 2026.

\bibitem{Kobayashi1987}
S.~Kobayashi, \emph{{Differential Geometry of Complex Vector Bundles}},
  Publications of the Mathematical Society of Japan, vol.~15, Iwanami Shoten
  and Princeton University Press, Tokyo and Princeton, NJ, 1987.

\bibitem{LeeSidiropoulos2010}
J.~R. Lee and A.~Sidiropoulos, \emph{{Genus and the geometry of the cut
  graph}}, {Proceedings of the Twenty-First Annual ACM--SIAM Symposium on
  Discrete Algorithms}, SIAM, 2010, pp.~193--201.

\bibitem{Liokumovich2014}
Y.~Liokumovich, \emph{{Surfaces of small diameter with large width}}, J. Topol.
  Anal. \textbf{6} (2014), no.~3, 383--396.

\bibitem{LiuPetri2026}
Mingkun Liu and Bram Petri, \emph{Random surfaces with large systoles},
  Mathematische Annalen \textbf{396} (2026), no.~1, 8.

\bibitem{LS}
A.~Lubotzky and Y.~Shalom, \emph{{Finite representations in the unitary dual
  and Ramanujan groups}}, Contemp. Math. \textbf{347} (2004), 173--189.

\bibitem{Maz24}
L.~Mazet, \emph{{Stable minimal hypersurfaces in $\mathbb{R}^6$}}, preprint,
  arXiv:2405.14676, 2024.

\bibitem{Micallef1984}
M.~J. Micallef, \emph{{Stable minimal surfaces in Euclidean space}}, J.
  Differential Geom. \textbf{19} (1984), no.~1, 57--84.

\bibitem{Micallef1986}
\bysame, \emph{{Stable minimal surfaces in flat tori}}, {Complex Differential
  Geometry and Nonlinear Differential Equations} (Brunswick, Maine, 1984),
  Contemp. Math., vol.~49, Amer. Math. Soc., Providence, RI, 1986, pp.~73--78.

\bibitem{MM1988}
M.~J. Micallef and J.~D. Moore, \emph{{Minimal two-spheres and the topology of
  manifolds with positive curvature on totally isotropic two-planes}}, Ann. of
  Math. (2) \textbf{127} (1988), no.~1, 199--227.

\bibitem{NS}
M.~S. Narasimhan and C.~S. Seshadri, \emph{{Stable and unitary vector bundles
  on a compact Riemann surface}}, Ann. of Math. (2) \textbf{82} (1965),
  540--567.

\bibitem{Pog81}
A.~V. Pogorelov, \emph{{On the stability of minimal surfaces}}, Dokl. Akad.
  Nauk SSSR \textbf{260} (1981), no.~2, 293--295, English transl. in
  \emph{Soviet Math. Dokl.} \textbf{24} (1981), 274--276.

\bibitem{Ramanan2009}
S.~Ramanan, \emph{{The theory of vector bundles on algebraic curves with some
  applications}}, {Moduli Spaces and Vector Bundles}, London Math. Soc. Lecture
  Note Ser., vol. 359, Cambridge University Press, Cambridge, 2009,
  pp.~165--210.

\bibitem{SagmanThalmaier2026}
N.~Sagman and T.-R. Thalmaier, \emph{{Holomorphicity of stable minimal surfaces
  of low genus}}, arXiv:2605.04399, 2026.

\bibitem{Schoen2006}
R.~Schoen, \emph{{Minimal submanifolds in higher codimension}}, Mat. Contemp.
  \textbf{30} (2006), 169--199.

\end{thebibliography}


\begin{thebibliography}{99}
\bibitem{AS}
L.~V.~Ahlfors and L.~Sario,
\emph{Riemann Surfaces},
Princeton Mathematical Series, vol.~26,
Princeton University Press, Princeton, NJ, 1960.

\bibitem{ArezzoMicallef2000}
C.~Arezzo and M.~J.~Micallef, \emph{Minimal surfaces in flat tori},
Geom. Funct. Anal. \textbf{10} (2000), 679--701.

\bibitem{AMP2002}
C.~Arezzo, M.~J.~Micallef, and G.~P.~Pirola,
\emph{Stable minimal surfaces of finite total curvature},
Comm. Anal. Geom. \textbf{10} (2002), no.~1, 11--22.

\bibitem{CKM2025}
D.~R.~Cheng, S.~Karigiannis, and J.~Madnick,
\emph{Holomorphicity of parabolic stable minimal surfaces of high
codimension}, arXiv:2509.22155, 2025, revised 2026.

\bibitem{ChernOsserman1967}
S.-S.~Chern and R.~Osserman,
\emph{Complete minimal surfaces in Euclidean $n$-space},
J. Analyse Math. \textbf{19} (1967), 15--34.


\bibitem{CL24} O.~Chodosh and C.~Li, Stable minimal hypersurfaces in $\mathbf{R}^4$, \emph{Acta Math.} \textbf{233} (2024), no.~1, 1--31.

\bibitem{CL23} O.~Chodosh and C.~Li, Stable anisotropic minimal hypersurfaces in $\mathbf{R}^4$, \emph{Forum Math. Pi} \textbf{11} (2023), Paper No.~e3, 22~pp.

\bibitem{CLMS26} O.~Chodosh, C.~Li, P.~Minter, and D.~Stryker, Stable minimal hypersurfaces in $\mathbf{R}^5$, \emph{Ann. of Math. (2)} \textbf{204} (2026), no.~1, 423--453.


\bibitem{CMR24} G.~Catino, P.~Mastrolia, and A.~Roncoroni, Two rigidity results for stable minimal hypersurfaces, \emph{Geom. Funct. Anal.} \textbf{34} (2024), no.~1, 1--18.



\bibitem{Demailly2012}
J.-P.~Demailly,
\emph{Complex Analytic and Differential Geometry},
Institut Fourier, Universit\'e de Grenoble~I,
online manuscript, version of June~21, 2012.


\bibitem{dCP79} M.~do Carmo and C.~K.~Peng, Stable complete minimal surfaces in $\mathbf{R}^3$ are planes, \emph{Bull. Amer. Math. Soc. (N.S.)} \textbf{1} (1979), no.~6, 903--906.

\bibitem{FCS80} D.~Fischer-Colbrie and R.~Schoen, The structure of complete stable minimal surfaces in $3$-manifolds of nonnegative scalar curvature, \emph{Comm. Pure Appl. Math.} \textbf{33} (1980), no.~2, 199--211.




\bibitem{Fraser2003}
A.~M.~Fraser, \emph{Fundamental groups of manifolds with positive
isotropic curvature}, Ann. of Math. (2) \textbf{158} (2003), no.~1,
345--354.

\bibitem{Fraser2026}
A.~Fraser, \emph{Stability theory for minimal surfaces},
in \emph{Proceedings of the International Congress of Mathematicians
2026}, vol.~4, SIAM, 2026, pp.~139--159.


\bibitem{FraserSchoen2025}
A.~Fraser and R.~Schoen, \emph{Stability and largeness properties of
minimal surfaces in higher codimension}, J. Eur. Math. Soc. (2025)

\bibitem{FraserWolfson2006}
A.~Fraser and J.~Wolfson,
\emph{The fundamental group of manifolds of positive isotropic
curvature and surface groups},
Duke Math. J. \textbf{133} (2006), no.~2, 325--334.


\bibitem{HLW26} H.~Hong, H.~Li, and G.~Wang, The stable Bernstein theorem in $\mathbb{R}^7$, preprint, 2026, arXiv:2609.15720.

\bibitem{HS85} R.~Hardt and L.~Simon, Area minimizing hypersurfaces with isolated singularities, \emph{J. Reine Angew. Math.} \textbf{362} (1985), 102--129.

\bibitem{Kobayashi1987}
S.~Kobayashi,
\emph{Differential Geometry of Complex Vector Bundles},
Publications of the Mathematical Society of Japan, vol.~15,
Iwanami Shoten, Tokyo, and Princeton University Press,
Princeton, NJ, 1987.

\bibitem{LeeSidiropoulos2010}
J.~R.~Lee and A.~Sidiropoulos,
\emph{Genus and the geometry of the cut graph},
in \emph{Proceedings of the Twenty-First Annual ACM--SIAM Symposium
on Discrete Algorithms}, SIAM, 2010, pp.~193--201.

\bibitem{Liokumovich2014}
Y. Liokumovich,
\emph{Surfaces of small diameter with large width},
J. Topol. Anal. \textbf{6} (2014), no.~3,
383--396.

\bibitem{LS}
A.~Lubotzky and Y.~Shalom,
\emph{Finite representations in the unitary dual and Ramanujan groups},
Contemp. Math. \textbf{347} (2004), 173--189.

\bibitem{Maz24} L.~Mazet, Stable minimal hypersurfaces in $\mathbb{R}^6$, preprint, 2024, arXiv:2405.14676.


\bibitem{Micallef1984}
M.~J.~Micallef, \emph{Stable minimal surfaces in Euclidean space},
J. Differential Geom. \textbf{19} (1984), no.~1, 57--84.

\bibitem{Micallef1986}
M.~J.~Micallef, \emph{Stable minimal surfaces in flat tori},
in \emph{Complex Differential Geometry and Nonlinear Differential
Equations} (Brunswick, Maine, 1984), Contemp. Math., vol.~49,
Amer. Math. Soc., Providence, RI, 1986, pp.~73--78.

\bibitem{MM1988}
M.~J.~Micallef and J.~D.~Moore, \emph{Minimal two-spheres and the topology
of manifolds with positive curvature on totally isotropic two-planes},
Ann. of Math. (2) \textbf{127} (1988), no.~1, 199--227.

\bibitem{NS}
M.~S.~Narasimhan and C.~S.~Seshadri,
\emph{Stable and unitary vector bundles on a compact Riemann surface},
Ann. of Math. (2) \textbf{82} (1965), 540--567.


\bibitem{Pog81} A.~V.~Pogorelov, On the stability of minimal surfaces, \emph{Dokl. Akad. Nauk SSSR} \textbf{260} (1981), no.~2, 293--295; English transl. in \emph{Soviet Math. Dokl.} \textbf{24} (1981), 274--276.

\bibitem{Qiu03} W.~Qiu, Stable minimal surfaces in higher dimensions, \emph{J. Differential Geom.} \textbf{63} (2003), no.~1, 1--18.

\bibitem{Ramanan2009}
S.~Ramanan,
The theory of vector bundles on algebraic curves with some applications,
in \emph{Moduli Spaces and Vector Bundles},
London Math. Soc. Lecture Note Ser., vol.~359,
Cambridge University Press, Cambridge, 2009, pp.~165--210.

\bibitem{SagmanThalmaier2026}
N.~Sagman and T.-R.~Thalmaier,
\emph{Holomorphicity of stable minimal surfaces of low genus},
arXiv:2605.04399, 2026.

\bibitem{Schoen2006}
R.~Schoen, \emph{Minimal submanifolds in higher codimension},
Mat. Contemp. \textbf{30} (2006), 169--199.
\end{thebibliography}

\end{document}